\documentclass{amsart}
\usepackage{amssymb,amsmath,amscd}
\newtheorem{theorem}{Theorem}[section]
\newtheorem{lemma}[theorem]{Lemma}
\newtheorem{corollary}[theorem]{Corollary}
\newtheorem{definition}[theorem]{Definition}

\newtheorem{remark}[theorem]{\it Remark}
\newtheorem{example}[theorem]{Example}
\newtheorem{proposition}[theorem]{Proposition}
\newtheorem{conjecture}[theorem]{Conjecture}
\newtheorem{problem}[theorem]{Problem}

\begin{document}


\def\sMr{\widetilde{M}_\mathbf{r}(\mathbb{C}P^m)}
\def\sMs{\widetilde{M}_{\Lambda-\mathbf{r}}(\mathbb{C}P^{n-m-2})}
\def\Mr{M_\mathbf{r}(\mathbb{C}P^m)}
\def\Nr{N_\mathbf{r}(\mathbb{C}P^m)}
\def\Ms{M_{\Lambda - \mathbf{r}}(\mathbb{C}P^{n-m-2})}
\def\Ns{N_{\Lambda - \mathbf{r}}(\mathbb{C}P^{n-m-2})}
\def\Grk{Gr_k(\mathbb{C}^n)}
\def\b{\mathbb}
\def\C{\mathbb{C}}
\def\R{\mathbb{R}}
\def\Z{\mathbb{Z}}
\def\br{\mathbf{r}}
\def\bs{\mathbf{s}}
\def\bk{\mathbf{k}}
\def\bl{\mathbf{l}}
\def\t{\mathrm}
\def\ba{\mathbf{a}}
\def\bb{\mathbf{b}}
\def\Lf{\mathcal{L}_\bk^\ba}
\def\Lg{\mathcal{L}_\bl^\bb}
\def\Fk{F_{\bk}(\C^n)}
\def\Fl{F_{\bl}(\C^n)}
\def\Rk{R_{\bk}}
\def\Rl{R_{\bl}}
\def\L{\mathbf{\Lambda}}
\def\Zo{\Z_{\scriptscriptstyle{\geq 0}}}
\def\Mu{\mathcal{M}_\br(\mathbb{CP}^m)}

\title[Duality of weight varieties]{The Chevalley Involution and a Duality of 
Weight Varieties}

\author{Benjamin J. Howard* and John J. Millson*}
\thanks{*Both authors were partially supported by NSF grant DMS-0104006}
\dedicatory{To the memory of Armand Borel}
\date{\today}

\begin{abstract}
In this paper we show that the classical notion of association
of projective point sets,\ \cite{DO}, Chapter III, is a special
case of a general duality between weight varieties (i.e.
torus quotients of flag manifolds) of a reductive group $G$ induced by the 
action
of the Chevalley involution on the set of these quotients. We compute 
the dualities explicitly on both the classical and quantum levels
for the case of  the weight varieties associated to $GL_n(\C)$. 
In particular we obtain the
following formula as a special case. Let $\br=(r_1,\ldots,r_n)$
be an $n$-tuple of positive real numbers and $M_{\br}(\mathbb{CP}^{m})$
be the moduli space of semistable weighted (by $\br$) configurations
of $n$ points in $\mathbb{CP}^{m}$ modulo projective equivalence, see
for example \cite{FM}. 
Let $\mathbf{\Lambda}$ be the vector in 
$\mathbb{R}^n$ with all components equal to $\sum_i r_i/(m+1)$. Then 
$M_{\br}(\mathbb{CP}^{m}) 
\cong M_{\mathbf{\Lambda} -\br}(\mathbb{CP}^{n-m-2})$
(the meaning of $\cong$ depends on $\br$ and will be explained below,
see Theorem \ref{association}).
We conclude by studying ``self-duality'' i.e. those cases where the duality
isomorphism carries the torus quotient into itself. We characterize when
such a self-duality is trivial, i.e. equal to the identity map. In
particular we show that all self-dualities are nontrival for
the weight varieties associated to the exceptional groups. The quantum 
version of this problem, i.e. determining for which self-adjoint 
representations $V$ of $G$ the Chevalley involution acts as a scalar
on the zero weight space $V[0]$,
is important in connection with the irreducibility of the
representations
of Artin groups of Lie type which are obtained as the monodromy of
the Casimir connection, see \cite{MillsonToledano}, and will be treated
in \cite{HowardMillsonToledano}.

\end{abstract}

\maketitle
\section{Introduction}

In this paper we show that the classical notion of association
of projective point sets,\cite{DO}, Chapter III, is a special
case of a general duality between weight varieties (i.e
torus quotients of flag manifolds) of a reductive group $G$ induced by the 
action
of the Chevalley involution on the set of these quotients. 
We will build up the theory in stages, first the duality for
Grassmannians of $GL_n(\mathbb{C})$, then for general flag manifolds
of $GL_n(\mathbb{C})$ then the duality for general flag manifolds
of general semisimple complex groups. At each stage there are three
types of isomorphism theorems, the first type is a K\"ahler isomorphism of 
symplectic
quotients, the second type is an algebraic isomorphism of Mumford quotients
and the third type is an explicit formula for the isomorphism of homogeneous
coordinate rings in terms of combinatorial Lie theory.  We now give details.

\subsection{Duality results for torus quotients of Grassmannians}

\subsubsection{Duality of symplectic quotients of Grassmannians}
We describe the duality result for symplectic quotients of Grassmannians.
Let $B$ be the nondegenerate symmetric bilinear form on
$\C^n$ such that the standard basis
$\{\epsilon_1,...,\epsilon_n \}$ is orthonormal.
Let $T$ be the maximal compact subgroup of $H$ (so $T$ is a product of
$\dim_\C(H)$ circles).
Then the
operation $\Psi$ of taking orthogonal complement with respect
to $B$ induces a map $\Psi: Gr_k(\C^n) \to Gr_{n-k}(\C^n)$
which carries torus orbits to torus orbits, since it satisfies the
formula 

\begin{equation} \label{torusequivariance} 
 \Psi(h \cdot x) = h^{-1} \cdot \Psi(x)
\end{equation}

Let $|\br| = \sum_{i=1}^n r_i$ and $a = |\br|/k$.   
We will see in \S \ref{Chevalleyinvolution}
that $\Psi$ maps the torus momentum level $\br =
(r_1,...,r_n)$ for $T$ to the torus momentum level $\L - \br
= (a - r_1,...,a -r_n)$.  Consequently we obtain
the following duality theorem.

\vskip 12pt

\begin{theorem}\label{symplecticquotients}
Let $\br = (r_1,...,r_n) \in \mathfrak{t}^* \cong \R^n$ 
be in the image of the momentum mapping
for the action of $T$ on $(Gr_k(\C^n), a \omega_k)$. Then the map
$\Psi:Gr_k(\C^n) \to  Gr_{n-k}(\C^n)$ induces a homeomorphism
(also to be denoted $\Psi$) of the symplectic quotients
$(Gr_k(\C^n), a\omega_k)//_\br T$ to 
$(Gr_{n-k}(\C^n), a \omega_{n-k})//_{\L -\br} T$. In
case $\br$ is not on a wall, see \cite{FM}, \S 4, then both symplectic quotients
are smooth and the map $\Psi$ is a K\"ahler isomorphism.

\end{theorem}

\subsubsection{Duality of Mumford quotients of Grassmannians}

\paragraph{\bf Linearization of the torus action}

Let $H$ be the maximal torus of $GL_n(\b{C})$ consisting of the
nonzero diagonal matrices and let $PH$ be the image of $H$ in
$PGL_n(\b{C})$. Let $\omega_k$ be the symplectic form on
$Gr_k(\b{C}^n)$ induced by the
embedding into $\mathfrak{u}(n)^*$ as the orbit of the $k$--th
fundamental weight $\varpi_k$,  the highest weight of the
representation of $GL_n(\C)$ on $\bigwedge^k(\b{C}^n)$. Let
$\mathcal{L}_k$ be the dual  of the $k$--th exterior power $\mathcal{T}_k$ of 
the
tautological $k$-plane bundle over $Gr_k(\b{C}^n)$. We will refer to
$\mathcal{T}_k$ as the tautological line bundle over $Gr_k(\b{C}^n)$.
Let $P \subset GL_n(\mathbb{C})$ be 
the stabilizer
of the coordinate plane spanned by the first $k$ standard basis vectors.
The bundle $\mathcal{L}_k$ is the homogeneous $ GL_n(\mathbb{C})$--bundle
over $Gr_k(\b{C}^n) = GL_n(\mathbb{C})/P$ with the isotropy
representation $det_k^{-1}:P \to \mathbb{C}^{\ast}$ where $det_k$ assigns 
to $p \in P$ the determinant of the upper $k$ by $k$ block of $p$.
In particular the total space of $\mathcal{L}_k$ is the quotient of the
product $GL_n(\mathbb{C}) \times \mathbb{C}$ by the equivalence relation
$(g,z) \sim (gp, det_k(p)z)$. 
Since $\mathcal{L}_k$ admits a $GL_n(\mathbb{C})$ action it admits an 
$H$ action. 
Since the entries of $\br$ are integers
we may identify $\br$ with the character $\chi_{\br}$ of $H$ 
whose value at the diagonal matrix with entries $(z_1,...,z_n)$ is
$z_1^{r_1}\cdots z_n^{r_n}$. For any integer $b$  we may use 
this character to twist the action of $H$ on the line bundle 
$\mathcal{L}_k^{\otimes b }$.   
We will use the symbol $\mathcal{L}_k^{\otimes b}(\br)$ to denote the $H$--line 
bundle $\mathcal{L}_k^{\otimes b}$ equipped with the twisted (by $\chi_{\br}$) 
$H$ 
action. The group $H$ acts on $Gr_k(\b{C}^n)$
through the quotient $PH$. In what follows we will need conditions on 
$\br$ that are necessary and sufficient in order that the action of $H$
on $\mathcal{L}_k^{\otimes b}(\br)$ descends to an action of $PH$. Let 
$|\br| = \sum_ir_i$. 

\begin{lemma}
The induced action
of $H$ on $\mathcal{L}_k^{\otimes b}(\br)$ descends to an action of
$PH$ if and  only if $bk = |\br|$. 
\end{lemma}
\begin{proof}
Let $h = \mu I$ be a nonzero scalar matrix. Then for $[g,z]$ in the
total space of $\mathcal{L}_k^{\otimes b}(\br)$ we have
$$h [g,z] = [hg, \chi_{\br}(h)z] = [gh, \mu^{|\br|} z] = 
[g, det_k(h)^{-b} \mu^{|\br|} z] = [g, \mu^{-kb + |\br|}z].$$ 
Thus $ h[g,z] = [g,z] \Leftrightarrow bk = |\br|$.
\end{proof}

The reader will verify that the condition $bk=|\br|$ is necessary in
order that there exist a nonzero section of $\mathcal{L}_k^{\otimes b}(\br)$
that is invariant under the group of nonzero scalar matrices. Thus
if it is not satisfied there will be no nonzero $H$--invariant sections
of $\mathcal{L}_k^{\otimes b}(\br)$. For this reason we assume that $|\br|$
is divisible by $k$ and we will reserve the symbol $a$ for the quotient
$|\br|/k$. We will abbreviate $\mathcal{L}_k^{\otimes a}$ to 
$\mathcal{L}_k^{a}$ henceforth. Let $\L = (a,a,...,a) \in \Z_+^n$.

\begin{definition}
For any integral $\br$ satisfying the condition that 
$|\br|$
is divisible by $k$ we will refer to the line bundle $\mathcal{L}_k^{a}(\br)$ 
equipped
with the previous action of $PH$ as the $\br$--linearization of the action of 
$PH$
on $Gr_k(\C^n)$.
\end{definition}

\paragraph{\bf The bundle isomorphism induced by the complex Hodge star}
Before stating our duality theorem concerning
Mumford quotients we need to recall the definition of a semistable point.
A point $x \in Gr_k(\C^n)$ is semistable for the
$\br$--linearization of $PH$ if there exists some $N \in \mathbb{Z}_+$
and an $H$--invariant section $s$ of $\mathcal{L}_k^{\otimes Na}(N\br)$
such that $s(x) \neq 0$ (here the symbol $\mathcal{L}(N\br)$ means we twist the
action of $H$ on $\mathcal{L}$ by the character $\chi_{\br}^N$).

Suppose now that $\br = (r_1,...,r_n) \in \mathbb{Z}_+^n$. 
We may form the Mumford quotients $Gr_k(\C^n)//_\br H$ and
$Gr_{n-k}(\C^n)//_{\Lambda -\br} H$ (using the linearizations
corresponding to $\br$ and $\Lambda-\br$ respectively). 
We will prove these Mumford quotients are isomorphic
varieties by constructing an explicit isomorphism of homogeneous 
coordinate rings. 
To this end choose a complex 
orientation of $\C^n$  and let
$\ast$ be the complex Hodge star (see \S 2)  associated to  
this orientation and the form $B$. We then define a bundle
isomorphism $\hat{\Psi}: \mathcal{L}_k \to \mathcal{L}_{n-k}$ as
follows. Let $x \in Gr_k(\C^n)$. Let $\tau \in \bigwedge^k((\C^n)^{\ast})$
and let $res_x:\bigwedge^k((\C^n)^{\ast}) \to \bigwedge^k(x^{\ast})$ be the
restriction map. Then we define $\hat{\Psi}$ by
$$\hat{\Psi}(res_x(\tau)) = res_{\Psi(x)}(\ast \tau).$$

It will be  proved in \S \ref{formulas} that $\hat{\Psi}$ is well defined 
and  it will be proved in \S \ref{ringtheorem} that 
 $\hat{\Psi}$ is an
$H$--morphism of line bundles (with the $H$ action on the target inverted).
The bundle map $\hat{\Psi}$ induces a map of sections 
$\widetilde{\Psi}:\Gamma(Gr_k(\C^n), \mathcal{L}_k) \to 
\Gamma(Gr_{n-k}(\C^n), \mathcal{L}_{n-k})$
given by
$$\widetilde{\Psi}(s)(x) = \hat{\Psi}(s(\Psi^{-1}(x))).$$
We will also  use $\hat{\Psi}$ resp. $\widetilde{\Psi}(x)$ to denote the maps 
on tensor powers that are induced by $\hat{\Psi}$ resp. $\widetilde{\Psi}$.
We will prove in  \S \ref{ringtheorem} that $\widetilde{\Psi}$ carries 
$H$--invariants to $H$--invariants and consequently induces an isomorphism
of homogeneous coordinate rings. We summarize these statements in

\begin{theorem}\label{rings}
\hfill
\begin{enumerate}
\item For all positive integers $N$ the   bundle isomorphism $\hat{\Psi}$ from 
$\mathcal{L}_k^{\otimes N a}(N \br)$ to 
$\mathcal{L}_{n-k}^{\otimes Na}(N (\L - \br))$
satisfies 
$$\hat{\Psi} \circ h = h^{-1} \circ \hat{\Psi}, h \in H.$$
\item The induced  isomorphisms of sections $\widetilde{\Psi}$ carry 
$H$--invariant 
sections of $\mathcal{L}_k^{\otimes Na}(N \br)$ to $H$--invariant 
sections of $\mathcal{L}_{n-k}^{\otimes N a} (N(\L - \br))$ for all $N$,
and consequently induce an isomorphism again to be denoted 
$\widetilde{\Psi}$ of the homogeneous coordinate rings of 
the Mumford quotients $Gr_k(\C^n)//_\br H$ and
$Gr_{n-k}(\C^n)//_{\L -\br} H$.
\end{enumerate}
\end{theorem}
  
 \vskip 12pt
 
\begin{remark}
We could have avoided the choice of form $B$ by defining
the duality map to be the map from the $k$--planes in $\C^n$
to the $(n-k)$--planes in $(\C^n)^*$ given by mapping a plane $x$
to its annihilator as suggested by \cite{Dolgachev}, Exercise 12.7, page 203.
However then we would not have had the explicit bundle map $\widetilde{\Psi}$
induced by the complex Hodge star.
\end{remark}

We will not give an explicit formula for the isomorphism $\widetilde{\Psi}$
on the standard basis for the homogeneous coordinates because it
is not much simpler than the more general formula of Theorem 
\ref{explicitformula}. Moreover the formulas in this case may be 
found in \cite{DO}, Ch. III.

\subsection{Duality of weighted projective configurations}
We now explain how the duality theorem for weighted projective configurations
(and its special case , the association of projective point sets) follows 
from the duality of symplectic torus quotients of Grassmannians.

\begin{theorem}\label{association}
 
Let $M_{\br}(\mathbb{CP}^m)$ be the moduli space of semistable weighted
(by $\br$) configurations of $n$--points on $\mathbb{CP}^m$. Then
 $$M_{\br}(\mathbb{CP}^m) \cong 
 M_{\mathbf{\Lambda} -\br}(\mathbb{CP}^{n-m-2}).$$ 
 By the symbol $\cong$ we mean the two spaces are homeomorphic   
and in case $\br$ is not on a wall, ( \cite{FM}, \S 4),then they are 
isomorphic as 
 K\"ahler manifolds. If $\br$ is integral then by $\cong$ we mean
 isomorphism as algebraic varieties.
 \end{theorem}
 
 \begin{proof} In what follows the symbol $\cong$ will have the same
 meaning as in the statement of the Theorem.
 By Gelfand-MacPherson duality, see for example \cite{FM}, \S 8, we have 
 $M_{\br}(\mathbb{CP}^m) \cong
 Gr_{m+1}(\C^n)//_{\br} H$ and $M_{\mathbf{\Lambda} -\br}(\mathbb{CP}^{n-m-2}) 
 \cong
 Gr_{n-m-1}(\C^n)//_{\mathbf{\Lambda}-\br} H$. But by Theorem 
 \ref{symplecticquotients} we have 
$Gr_{m+1}(\C^n)//_{\br} H \cong Gr_{n-m-1}(\C^n)//_{\mathbf{\Lambda}-\br} H$.
\end{proof}

\begin{remark} If all the components of $\br$ are equal (the 
``democratic linearization'') then the resulting isomorphism is
the classical association isomorphism, see \cite{DO}, Ch. III.
\end{remark}

\subsection{Duality results for torus quotients of flag manifolds}

It is remarkable that the  duality theorems for Mumford 
quotients and
symplectic quotients for Grassmannians almost immediately imply the
corresponding results for torus quotients of general flag manifolds.
The proofs are based on using the following diagram to promote
the duality theorem for Grassmannians to a duality theorem for
flag manifolds. We let $\Psi:F_\bk(\C^n) \to F_\bl(\C^n)$ be the
mapping that makes the following diagram commute.

\begin{center}
\begin{math}
\begin{CD}
F_\bk(\C^n) @>{\Psi}>> F_\bl(\C^n) \\
@V{i}VV                  @VV{i}V \\
 \prod_{i \leq m} Gr_{k_i}(\C^n) 
@>F \circ {\prod_i \Psi_i}>> 
\prod_{i \leq m} Gr_{l_i}(\C^n) 
\end{CD}
\end{math}
\end{center} 

\bigskip

Here $\bk=(k_1,...,k_m)$ resp. $\bl=(l_1,...,l_m):= (n-k_m,...,n-k_1)$ and $\Fk$ 
resp. $\Fl$
denotes the manifold of flags consisting of subspaces of dimensions
$k_1< \cdots k_m$ resp. $l_1 < \cdots < l_m$ and
$F$ is the map on products of Grassmannians that reverses the order of the 
factors. We assume $\Fk$ is given the
symplectic structure inherited by embedding it as the coadjoint orbit
of $\lambda = a_1\varpi_{k_1}+ \cdots a_m \varpi_{k_m}$. We put $\ba = (a_1,...,a_m)$,
$\bb=(b_1,...,b_m):=(a_m,...,a_1)$ and $|\ba|=\sum a_i$.
We let $\mathbf{\Lambda} = (|\ba|,\cdots,|\ba|)$ and
$\bs = \mathbf{\Lambda} -\br$.

\subsubsection{The duality theorem for symplectic quotients of flag manifolds.}

First the duality theorem for symplectic quotients (the proof may be
found in \S 2).

\begin{theorem}
The map $\Psi$ induces a homeomorphism of  symplectic quotients:
$$\overline{\Psi}:  F_\bk(\C^n) //_\br T \to F_\bl(\C^n)//_\bs T.$$
Furthermore, if $\br$ is a regular value of the momentum mapping, then
the symplectic quotients are K\"ahler manifolds and $\overline{\Psi}$ is a
K\"ahler isomorphism.  
\end{theorem}

Now the duality theorem for Mumford quotients. Let $\Lf$ be the homogeneous
line bundle with isotropy representation the character corresponding to
the negative of the dominant weight $\lambda = a_1\varpi_{k_1}+ \cdots + a_m 
\varpi_{k_m}$.
We let $\br \in \mathbb{Z}^n_{+}$ and let $\Lf(\br)$ be the $H$-bundle
represented by the line bundle $\Lf$ with the $H$ action  twisted by
the character corresponding to $\br$.
We first give the relation between $\ba$ and $\br$ that is necessary and 
sufficient in order that $PH$ act on $\Lf(\br)$. The following lemma
is Lemma \ref{existenceoflinearization} in the text. 

\begin{lemma}
The action of $H$ on $\Lf(\br)$ descends to an action of $PH$ iff
$|\br| = \sum_i a_i k_i$.
\end{lemma}

We define the bundle map $\hat{\Psi}$ to be the tensor product
of the bundle maps for the Grassmannians followed by a reversal
of tensor factors.

We now state our isomorphism theorem for Mumford quotients, the proof
is to be found in \ \S 4.

\begin{theorem}
\begin{enumerate}
\item The map $\widetilde{\Psi}$ induces an isomorphism of graded rings:
$$ \bigoplus_{N=0}^\infty 
\Gamma(F_\bk(\C^n),\mathcal{L}_\bk^{N \ba}(N \br))^H \cong 
\bigoplus_{N=0}^\infty 
\Gamma(F_\bl(\C^n),\mathcal{L}_\bl^{N \bb}(N \bs))^H.$$
\item Equivalently, the  map $\widetilde{\Psi}$ induces an isomorphism of 
Mumford quotients: 
$$F_\bk(\C^n)//_\br H \cong F_\bl(\C^n)//_\bs H.$$
\end{enumerate}
\end{theorem}

\paragraph{\bf An explicit formula for the ring isomorphism $\widetilde{\Psi}$}

 We next give an explicit formula for $\widetilde{\Psi}$ on the
 homogeneous coordinate ring of the flag manifold in terms
 of semistandard Young tableaux. The proof of  the following theorem
 is to be found in \S 5.

\vskip 12pt

Let $\lambda = \sum_i a_i \varpi_k$.  
The $N$-th graded summand $R_\bk^{(N)}$ of the homogeneous coordinate ring of 
$F_\bk(\C^n)$ is given by
$$R_\bk^{(N)} = 
\Gamma(F_\bk(\C^n), (\Lf)^{\otimes N})
= V_{N \lambda}.$$
Furthermore the $N$--th graded summand $(R_\bk^{(N)})^H$ of the homogeneous 
coordinate ring of the Mumford quotient $F_\bk(\mathbb{C}^n)//_{\br} H$ is 
given by
$$(R_\bk^{(N)})^H = 
\Gamma(F_\bk(\C^n), (\Lf)^{\otimes N}(N\br))^H
= V_{N \lambda}(N\br).$$

The last symbol denotes the $N \br$--th weight space of the irreducible 
representation of $GL(n,\mathbb{C})$ with highest weight $N \lambda$.
It is a standard result
in representation theory \cite{Boerner}, Ch. V, Theorem 5.3, that there is a 
basis for $V_{N \lambda}$ resp.
$V_{N \lambda}(N \br)$ consisting of the semistandard fillings resp.
semistandard fillings 
of weight $N\br$ of the Young diagram $D_\ba$ (the $i^{th}$ row has
length equal to the $i^{th}$ component of $N \lambda$)
by the integers between $1$ and $n$ inclusive.  We will call this the standard
basis of $R_\bk^{(N)})$ (resp. $(R_\bk^{(N)})^H$).  If $T$ is a semistandard 
filling (of any weight) 
of the above Young diagram by $1,2,..,n$ we will
let $f_T$ denote the corresponding element of the homogeneous coordinate ring
of the flag manifold $F_\bk(\C^n)$. The set of all such $f_T$ is a
basis for $R_\bk^{(N)}$. Note that $\widetilde{\Psi}$ induces a map from
$R_\bk^{(N)}$ to $R_\bl^{(N)} = \Gamma(F_\bl(\C^n), (\Lg)^{\otimes N})$.

We will describe this map relative to the
standard basis, a fortiori this will describe the map on the subrings of
$H$--invariants.

We now describe a map from semistandard tableaux of weight $\br$ on the diagram
$D_\ba$ to semistandard tableaux of weight 
$\L - \br$ on the diagram $D_\bb$
which we will denote $T \mapsto \ast T$.  We explain how
to obtain the dual tableau $\ast T$ with an example.

\begin{example}
Let $\lambda = 2 \varpi_2 + \varpi_3, N = 1$.
\begin{center}
$T=$
\begin{tabular}{| c | c | c |}
\hline
2  & 1  & 2  \\ \hline
3  & 4  & 5  \\ \hline
5  \\ \cline{1-1}
\end{tabular}
$\Longrightarrow \widetilde{T} = $
\begin{tabular}{| c | c | c |}
\hline
    2  &     1  &     2  \\ \hline
    3  &     4  &     5  \\ \hline
    5  & \it{2} & \it{1} \\ \hline
\it{1} & \it{3} & \it{3} \\ \hline
\it{4} & \it{5} & \it{4} \\ \hline
\end{tabular}
$\Longrightarrow  *T = $
\begin{tabular}{| c | c | c |}
\hline
1  &  2  &  1   \\ \hline
3  &  3  &  4   \\ \hline
4  &  5  \\ \cline{1-2}
\end{tabular}

\end{center}
In general, extend $T$ to a rectangular $n$ by $|\ba|$ diagram and fill in the 
complementary indices in each column, listed in increasing order.  
Write the added columns in reverse order to get $\ast T$.
\end{example}

\vskip 12pt
We attach a sign $\epsilon_T$ to each semistandard tableau as follows. Form
the enlarged tableau $\widetilde{T}$ as above.  For each column $C_i$
of $\widetilde{T}$ define $\epsilon_i$ to be the sign of the permutation
of $1,2,...,n$ represented by that column read from top to bottom.
Then define 
$$\epsilon_T = \epsilon_1...\epsilon_n.$$

\begin{theorem}\label{tableaux}\label{explicitformula}
 The $N$-th graded component of the isomorphism $\widetilde{\Psi}$
is diagonal relative to the standard bases for $R_\bk^{(N)}$ and 
$R_\bl^{(N)}$.   

Moreover we have the formula
$$\widetilde{\Psi}(f_T) = \epsilon_T  f_{\ast T}.$$
\end{theorem}

\begin{remark} It is not obvious that the map of graded vector spaces 
given by the formula in the theorem defines a {\emph ring} homomorphism.
Even for the case of rectangular Young diagrams and when the weights 
$r_i$ are all equal, 
a direct algebraic proof using the 
Pl\"ucker relations is not easy and was given in the thesis of D.\ Ortland - 
see the proof of Chapter III,Theorem 1, in \cite{DO}. However we know
a fortiori that this map of semistandard tableaux is induced by
the ring isomorphism $\widetilde{\Psi}$.
\end{remark}

\subsection{Duality for general semisimple complex Lie groups}
Let $\theta$ be the Chevalley involution of the Lie group $GL_n(\mathbb{C})$
whence 
$$\theta(g) = (g^{t})^{-1}.$$
Then $\theta$ carries a standard parabolic subgroup $P$ to its opposite
$P^{opp}$ and induces a map $\Theta:G/P \to G/P^{opp}$ given by

\begin{equation} \label{dualitymap}
\Theta(gP) = \theta(g) P^{opp}.
\end{equation}

Next if $\chi$ is a character of $P$ then the 
$\chi^{\theta} : = \chi \circ \theta$ is a character of $P^{opp}$.
Let $\mathcal{L}_{\chi}$ and $\mathcal{L}_{\chi^{\theta}}$ be the corresponding
homogeneous line bundles. Then we obtain an isomorphism of line bundles
$\widehat{\Theta}$ by defining
$$\widehat{\Theta}([g,z]) = [\theta(g),z].$$

We will prove the following lemma in \S 2, see 
Lemma \ref{relationtoChevalleyinvolution}.

\begin{lemma}
\hfill
\begin{enumerate}
\item $\Theta = \Psi.$
\item $\hat{\Theta} = \hat{\Psi}.$
\end{enumerate}
\end{lemma}

The critical point here is that with this formulation  i.e. using
Equation \ref{dualitymap}  extends the duality map
to a duality map of weight varieties for all reductive groups. To avoid
complications
we will restrict ourselves to the cases that either $G$ is semisimple
or $G = GL_n(\C)$ in what follows.

\begin{remark}\label{anotherdescription}
In fact we get the same map $\Theta$ (for general $G$) using another description.
Let $n(w_0) \in N(T)$ be a representative for the longest element
$w_0$ in the Weyl group. We may assume that it is fixed under
$\theta$ (see \S 2). Then we define $R:G/P^{opp} \to G/Q$
(where $Q$ is the standard parabolic conjugate to $P^{opp}$)
by $R(gP^{opp}) = g n(w_0) Q$. The reader will verify
that $R$ induces the identity on the flag manifold $M^{opp}$. 
Thus if we postcompose $\theta$ by $R$ we obtain the
same map $\Theta$ but we have another presentation in terms
of coset spaces. We will abuse notation and also use the same symbol $\Theta$
for this new presentation. We have then
$\Theta:G/P \to G/Q$ with 
$$\Theta(gP) = \theta(g) n(w_0) Q.$$
The reader will verify that with this description 
the induced map on line bundles carries $\mathcal{L}_{\chi_{\lambda}}$
to $\mathcal{L}_{\chi_{\lambda^{\vee}}}$ where $\lambda^{\vee}$
is the weight contragredient to $\lambda$. In what follows we will 
use whichever description of $\Theta$ is convenient.
\end{remark}

We now state two theorems and a conjecture, the analogues of the three theorems 
above
for the quotients of flag manifolds of the group $GL_n(\mathbb{C})$.

\subsubsection{Duality of symplectic quotients}
The following theorem is proved in \S 2.

\begin{theorem}
Let $K$ be a semisimple compact Lie group with complexification $G$ and $T$ be a 
 maximal
torus with $T \subset K$. Choose a Chevalley involution $\theta$ of $G$ such that $\theta$
carries $K$ into itself and satisfies $\theta(t) = t^{-1}$ for $t \in T$.
Let $S$ be a subtorus of $G$ and $Z(S)$ be the centralizer of $S$
whence $\theta(Z(S)) = Z(S)$. Let
$M=K/Z(S)$. Let $\br$
be an element of the moment polyhedron for the action of $T$ on $K/Z(S)$.
Then the Chevalley involution induces an isomorphism of K\"ahler manifolds
$$\overline{\Theta}: M//_{\br}T \to M//_{-\br} T.$$
\end{theorem}

\subsubsection{Duality of Mumford quotients}

We next state the corresponding general duality result for Mumford quotients.
We continue with the notation of the previous theorem. Let $H$ be the
complexification of $T$. There exists a Borel subgroup $B$ of $G$ such that
$B \cap \theta(B) = H$. We let $B^{opp}$ denote $\theta(B)$. We let $P$
be a parabolic subgroup of $G$ containing $B$ and let $P^{opp}= \theta(P)$. 
We will assume $G$ is the simply-connected group thus the character
lattice of $H$ is the weight lattice of the Lie algebra $\mathfrak{g}$. To emphasize
this point we make the definition

\begin{definition}
In what follows we say the symplectic manifold M is {\em integral} will mean 
$\lambda$ is in the weight lattice of $\mathfrak{g}$. Here $M$ corresponds
to the orbit of $\lambda$. We say $\br$ is {\em integral} if $\br$
is in the weight lattice.
\end{definition}
Let $\lambda$ be a dominant weight and assume that
the corresponding character $\chi_{\lambda}$ of  $H$ extends to a character of $P$
but does not extend to any larger parabolic. Let $\mathcal{L}_{\chi_{\lambda}}$
be the homogeneous line bundle over the flag manifold $M = G/P$ with isotropy 
representation
$\chi_{\lambda}^{-1}$. Then $\mathcal{L}_{\chi_{\lambda}}$ is a very ample $H$--bundle.
We let $\br$ be another weight and twist the action of $H$ on 
$\mathcal{L}_{\chi_{\lambda}}$ by the character of $H$ associated to $\br$
to obtain $\mathcal{L}_{\chi_{\lambda}}(\br)$.
We can then form the Mumford quotient $M//_{\br}H$ in the usual way. 
Similarly we obtain a very ample $H$--bundle 
$\mathcal{L}_{\chi_{\lambda \circ \theta}}(-\br)$ over 
the
flag manifold $M^{opp} = G/P^{opp} = G/Q$
with isotropy representation $\chi_{\lambda \circ \theta}$. We will denote the corresponding
Mumford quotient by $M^{opp}//_{-\br}H$. We have a bundle isomorphism
$\hat{\Theta}:\mathcal{L}_{\chi_{\lambda}}(\br) \to 
\mathcal{L}_{\chi_{\lambda \circ \theta}}(-\br)$ defined as
in the case of $GL_n(\mathbb{C})$. We obtain

\begin{theorem}
The bundle isomorphism $\hat{\Theta}$ induces an isomorphism of Mumford
quotients
$$M//_{\br}H \cong M^{opp}//_{-\br}H.$$
\end{theorem}

\subsubsection{An explicit formula for duality on the ring level}

There should be an explicit formula for computing the isomorphism
$\widetilde{\Theta}$ of homogeneneous coordinate rings associated to
the previous isomorphism in terms of the
Littelmann path model \cite{Littelmann2} (or any other model) for the irreducible representations
given by the graded summands of the two 
coordinate rings.

\begin{conjecture} \label{Littelmannpathreversal}
In the  Littelmann path models for the two homogeneous coordinate rings
the isomorphism $\widetilde{\Theta}$ is given by reversing  
Lakshmibai- Seshadri paths and translating the initial points of the reversed 
paths to the origin.

\end{conjecture}

\subsection{Self-duality}

We now suppose that the flag manifold 
$M$ and the level $\br$ have been chosen
so that $\overline{\Theta}$ carries $M//_{\br}H$ into itself. We will then say
the torus quotient is {\em self-dual}. We may then ask 

\begin{problem}[Classical Problem] For which self-dual torus quotients 
$M//_{\br} H$ is
$\overline{\Theta}$ is equal to the identity? 
\end{problem}

The quantum version of
the previous question is 

\begin{problem}[Quantum Problem]
For which  self-dual irreducible representations $V_{\lambda}$ does the
Chevalley involution act as a scalar on the zero weight space 
$V_{\lambda}[0]$.
\end{problem}

\begin{remark} If the Chevalley involution is inner then it is clear
that it acts on $V_{\lambda}[0]$, if not then the action on $V_{\lambda}$
and $V_{\lambda}[0]$ is defined only up to a scalar multiple, see 
\cite{MillsonToledano}, \S 4.3.
\end{remark}

The motivation for this problem is explained in \cite{MillsonToledano}
where it is solved for the groups $SL_n(\C)$ and $G_2(\C)$. For each
irreducible $V_{\lambda}$ as above the authors in 
\cite{MillsonToledano} construct an action of the
Artin group $B_{\mathfrak{g}}$ associated to  $\mathfrak{g}$ (the fundamental group
of the quotient by the Weyl group of the space of regular elements
in a Cartan subalgebra) on $V_{\lambda}[0]\otimes \C[[h]]$. This representation
is the monodromy representation of the Casimir connection and
by a very recent theorem of Toledano Laredo, see \cite{ToledanoLaredo2}, 
coincides with
the representation constructed by Lusztig, \cite{Lusztig}, Ch. 41, using the 
theory of 
quantum groups, see also \cite{ToledanoLaredo1} where the result was proved
for the case of $SL_n(\C)$. In \cite{MillsonToledano} the authors
began a study of the irreducibility of these representations.
The starting point of this study of irreducibility 
was the observation that in case $V_{\lambda}$ was self-dual these representations commute with
the action of the Chevalley involution and hence if the
Chevalley involution does not act as a scalar (as is nearly always
the case) they are reducible.

We now describe our solution of the classical problem.
The solution of the problem
for the groups $GL_n(\C)$ and $G_2$  follows from  earlier work of
the second author and V.\ Toledano Laredo, \cite{MillsonToledano}.
In what follows note that $\br = 0$ for all cases except for (1).

\begin{theorem}
\hfill

\begin{enumerate}
\item Suppose $G =GL_n(\mathbb{C})$. Assume that $\bk$ and $\br$ satisfy the 
self-duality conditions  $\bk = \bl$ and $\br = \bs$.
The  self-duality $\overline{\Theta} : 
F_\bk(\C^n) //_\br H \to F_\bk(\C^n)//_\br H$ is 
equal to the identity
if and only if  the flag manifold is 
\begin{enumerate}
\item $\mathbb{CP}^1$ with the symplectic form corresponding to $a\varpi_1$
and $\br = (a/2) \varpi_2$.
\item  $F_\bk(\C^n) = F_{1,n-1}(\mathbb{C}^n)$ with the
symplectic form $a \varpi_1 + a \varpi_{n-1}$ and $\br = a \varpi_n$.
\item $F_\bk(\C^n) = Gr_2(\mathbb{C}^4)$ with the
symplectic form $2a\varpi_2, a \in \mathbb{N}$ and $\br = a \varpi_4$.
\end{enumerate}
\item Suppose $G = Sp_{2n}(\C)$.Then
the duality map $\overline{\Theta}$ is  equal to the identity
if and only if the flag manifold is 
\begin{enumerate}
\item The projective space $\mathbb{CP}^{2n-1}$
with the symplectic form corresponding to a multiple of $\varpi_1$.
\item The Lagrangian Grassmannian $Gr_{2}^{0}(\C^4)$ with the symplectic
form corresponding to a multiple of $\varpi_2$.
\end{enumerate}

\item Suppose that $G= SO_{2n+1}(\C)$. Then the duality map 
$\overline{\Theta}$  is  equal to the identity
if and only if the flag manifold is  
\begin{enumerate}
\item The quadric hypersurface $\mathcal{Q}_{2n-1}\subset \mathbb{CP}^{2n}$
with the symplectic form corresponding to a multiple of $\varpi_1$.
\item The Lagrangian Grassmannian $Gr_{2}^{0}(\C^5)$
with the symplectic
form corresponding to a multiple of $\varpi_2$.
\end{enumerate}

\item Suppose that $G= SO_{2n}(\C)$. Then the duality map 
$\overline{\Theta}$ 
 is equal to the identity
if and only if the flag manifold is 
\begin{enumerate}
\item The quadric hypersurface $\mathcal{Q}_{2n-2}\subset \mathbb{CP}^{2n-1}$
with the symplectic form corresponding to a multiple of $\varpi_1$.
\item One of the  Lagrangian Grassmannians $Gr_{2}^{0}(\C^4)^{+} \cong \mathbb{CP}^1$
and $Gr_{2}^{0}(\C^4)^{-}\cong \mathbb{CP}^1$
with the K\"ahler forms corresponding to a constant curvature
form.
\item The Grassmannian of isotropic two-planes $Gr_2^0(\C^6)$.
\item One of the Lagrangian Grassmannians $Gr_{4}^{0}(\C^8)^{+}$ and 
$Gr_{4}^{0}(\C^8)^{-}$ with the symplectic forms corresponding to multiples
of $\varpi_3$ and $\varpi_4$.
\item The isotropic flag manifold $F_{1,2}^0(\C^{4}) 
\cong \mathbb{CP}^1 \times \mathbb{CP}^1$ with the K\"ahler
form corresponding to the sum of any two constant curvature forms. 
\end{enumerate}

\end{enumerate}

\end{theorem}

\begin{remark}
We see from the above that there are three infinite families of examples
where the duality is trivial, the line-hyperplane pairs, the lines
in symplectic vector spaces and the quadrics for the orthogonal groups.
It is remarkable that {\em all} the other examples are obtained from
these three infinite families using exceptional isomorphisms. We discuss
two examples in detail. First, the  example (4) (c), the isotropic
Grassmannian $Gr_2^0(\C^6)$ ,
is explained by the exceptional isomorphism
$D_3 \cong A_3$ which carries $Gr_2^0(\C^6)$ to the member of the
infinite family of line-hyperplane pairs given by $F_{1,3}(\C^4)$. 
Second, the example of the two Lagrangians in (4)(d), is explained by
triality. Indeed we have
$$Gr_4^0(\C^8)^+ \cong Gr_4^0(\C^8)^- \cong \mathcal{Q}_6.$$
As for the other two not-quite-obvious examples, (2)(b) and (3)(b), 
we have, using the exceptional isomorphism $C_2 \cong B_2$,
\begin{enumerate}
\item  \ $Gr_2^0(\C^4) \cong \mathcal{Q}_3$.
\item  \ $Gr_2^0(\C^5) \cong \mathbb{CP}^3$.
\end{enumerate}
\end{remark}

\subsection{The idea of the proof}

We conclude the mathematical part of this Introduction with a sketch
of the proof of the previous
theorem. The following definition
is critical in what follows.

\begin{definition} A representation $V_{\lambda}$ is {\em good}
if it is self-dual and for some $N$ the Chevalley involution $\theta$
does not act as a scalar on $V_{N\lambda}[0]$.
\end{definition}

We prove in what follows that if $M//_0 H$ is a self-dual torus quotient then
the duality map $\overline{\Theta}$ is nontrivial if and only if
the representation $V_{\lambda}$ is good, Theorem \ref{quantumtoclassical}.
Here we assume that $M$ is
the flag manifold associated to $\lambda$. We also prove that
the good representations are closed under Cartan products, see Definition
\ref{defCartanproduct}.
In fact any Cartan product of self-dual representations is good provided
at least one factor is good, 
Theorem \ref{Cartanproduct}. 
Thus,for example, to prove that all self-dualities are nontrivial for a given
group $G$ such that $-1 \in W$ we have only to prove that all
the fundamental representations are good. We prove this for
the groups $G_2, F_4, E_7$ and $E_8$ by branching
a fundamental representation to a carefully chosen maximal subgroup of maximal rank
and observing that this restriction contains either a good representation
or a nonself-dual representation. 

\bigskip
{\bf Acknowledgements}
We would like to thank Philip Foth who after reading an early version of
this paper showed us
a duality map different from our version of duality that led us to 
find the connection of our original duality with
the Chevalley involution. We would like to thank Jeffrey Adams, Shrawan Kumar and
Valerio Toledano Laredo for helpful conversations. Also we used
the program LiE to prove the fundamental representations of the exceptional
groups and the second fundamental representation of $SO(8)$ were good. 

The second author would like to acknowledge the fundamental role that
Armand Borel played in his career both as a mentor and as an example
of what it means to be a mathematician.

\section{The Chevalley involution and duality of symplectic quotients}
\label{Chevalleyinvolution}
We recall the definition of a Chevalley involution. Choose a Cartan
subalgebra $\mathfrak{h}$ and a Borel $\mathfrak{b}$ containing 
$\mathfrak{h}$. Thus we obtain a system of roots $R$ together with a positive
subsystem $R_+ \subset R$ and a simple subsystem $S \subset R_+$. For each
simple root $\alpha$ choose a root vector $x_{\alpha}$ corresponding to
$\alpha$.  Let $h_{\alpha} \in \mathfrak{h}$ be the coroot corresponding to
$\alpha$. Then there is a unique negative root vector $x_{-\alpha}$
such that $[x_{\alpha}, x_{-\alpha}] = h_{\alpha}$. We then have
the following consequence of the Chevalley presentation of $\mathfrak{g}$.

\begin{lemma}\label{Chevalley}
There exists a unique involutive automorphism $\theta$ of $\mathfrak{g}$
such that
 
\begin{enumerate}
\item $ \theta(x_{\alpha}) = - x_{-\alpha} \ \text{for all} \ \alpha \in S.$
\item $ \theta(x_{- \alpha}) = - x_{\alpha} \ \text{for all} \ \alpha \in S.$
\item $\theta(h_{\alpha}) = - h_{\alpha} \ \text{for all} \ \alpha \in S.$ 
\end{enumerate}

\end{lemma}

We will say a holomorphic involution of a simple complex Lie algebra
$\mathfrak{g}$ is a Chevalley involution if $\theta$ satisfies the
above formulas for some $\mathfrak{h}$, $\mathfrak{b}$ and vectors
$x_{\alpha}$,$x_{-\alpha}$ and $h_{\alpha}$, $\alpha \in S$.

\begin{remark}
Any two Chevalley involutions $\theta_1$ and $\theta_2$ are conjugate
by Proposition 2.8 of \cite{AdamsBarbaschVogan}. Furthermore it is
possible to choose root vectors $x_{\alpha}$ for all positive roots
$\alpha$ such that $(1)$ and $(2)$ continue to hold. This follows provided
one has chosen the structure constants $N_{\alpha,\beta}$  for 
$\mathfrak{g}$ in the Chevalley basis so that
$$N_{\alpha,\beta} = - N_{-\alpha,-\beta}.$$
See \cite{Samelson}, Ch II, \S 9.
\end{remark}

In our study of self-duality for the weight varieties associated to
the exceptional groups we will need the following lemma.

\begin{lemma}\label{characterizationofChevalley}
Suppose $H$ is a Cartan subgroup of a simple complex Lie group $G$. Suppose
$\theta_1$ and $\theta_2$ are holomorphic involutions of $G$ which
carry $H$ into itself and satisfy $\theta_i(h) = h^{-1}, h \in H, i=1,2$.
Then there exists $h \in H$ such that 
$$ \theta_2 = Ad h \circ \theta_1 \circ Ad h^{-1}.$$
Moreover both $\theta_1$ and $\theta_2$ are Chevalley involutions.
\end{lemma}
\begin{proof}
Let $\alpha$ be a root and $\mathfrak{g}_{\alpha}$ be the corresponding
root space. Then we have
$$\theta_i(\mathfrak{g}_{\alpha}) = \mathfrak{g}_{-\alpha}, i=1,2.$$
Let $\alpha_i,1 \leq i \leq l$ be the simple roots. For each
simple root $\alpha_i$ choose a Chevalley basis vector $x_{\alpha_i}$.
Since $\theta_1(x_{\alpha_i})$ and $\theta_2(x_{\alpha_i})$ both lie in
the one dimensional space  $\mathfrak{g}_{-\alpha}$ there exist
complex numbers $c_i$ such that 
$$\theta_2(x_{\alpha_i}) = c_i \theta_1(x_{\alpha_i}).$$
Choose $h \in H$ such that $Ad h(x_{\alpha_i}) = \sqrt{c_i}x_{\alpha_i}$
whence $Ad h(x_{-\alpha_i}) = (1/\sqrt{c_i}) x_{-\alpha_i}$.
Then $\theta_2 = Ad h^{-1} \circ \theta_1 \circ Ad h$.
\end{proof}

In what follows we will have a distinguished Cartan $H$. The above lemma
allows us to make an abuse of language and refer to {\em the}
Chevalley involution of $G$ (and $H$).

\subsubsection{The action of the Chevalley involution on a 
self-dual representation} \label{actionofChevalley}

In this subsection we recall
how $\theta$ acts on the weight space $V_{\lambda}[0]$ of a self-dual
representation $V_{\lambda}$, see
\cite{MillsonToledano}, \S 4.3. Indeed because $V_{\lambda}$ is self-dual
there exists $\Theta_{V_{\lambda}} \in Aut(V_{\lambda})$ of order $2$ which intertwines the action
$\rho$ of $GL_n(\C)$ with the action $\rho^{\theta} = \rho \circ \theta$
on the same space. By Schur's Lemma $\Theta_{V_{\lambda}}$ is unique up to multiplication
by $\pm 1$. Then the action of $\theta$ on $V_{\lambda}$ is defined to
be the action of $\Theta_{V_{\lambda}}$. It is then immediate that 
$\Theta_{V_{\lambda}}$ carries the zero weight space $V_{\lambda}[0]$
into itself. 

\begin{lemma}\label{agreementofChevalley}
Suppose $H$ is a Cartan subgroup of a simple complex Lie group $G$. Suppose
$\theta_1$ and $\theta_2$ are holomorphic involutions of $G$ which
carry $H$ into itself and satisfy $\theta_i(h) = h^{-1}, h \in H, i=1,2$
Let $V$ be a self-dual irreducible representation of $G$ and let
$\Theta^{(1)}_V$ and $\Theta^{(2)}_V$ be the operators assigned
to $\theta_1$ and $\theta_2$ according to the rule explained in
the preceding paragraph. Then $\Theta^{(1)}_V$ and $\Theta^{(2)}_V$
are conjugate by an element of $H$ acting on $V$ and consequently
the restrictions of $\Theta^{(1)}_V$ and $\Theta^{(2)}_V$ to
$V[0]$ coincide.
\end{lemma}

\subsubsection{The duality map}

Let $P$ be a standard parabolic subgroup of $G$ and $M$ be the 
flag manifold $G/P$. Let $P^{opp} = \theta(P)$ and 
$Q$ be the standard parabolic subgroup conjugated to $P^{opp}$.
Let $M^{opp} = G/P^{opp}$ and $N=G/Q$. Then $M^{opp} = N$.
We have defined the map $\Theta: M \to M^{opp}$  
$$\Theta(gP) = \theta(g)P^{opp}.$$
Equivalently we have defined the duality map $\Theta: M \to N$ by
$$\Theta(gP) = \theta(g) n(w_0)Q.$$

It is immediate that 
\begin{equation}\label{equivariance}
\Theta(gx) = \theta(g) \Theta(x), x \in M
\end{equation}

Let $v \in \mathfrak{g}$ and $V_M$ be the fundamental vector field
associated to $v$. Then the infinitesimal version of Equation 
(\ref{equivariance}) 
follows immediately from  Equation (\ref{equivariance}). We will need it below
so we state it as a lemma.

\begin{lemma}\label{infinitesimalequivariance}
$\Theta_{\ast}(V_M)$ is the fundamental vector field $\theta(v)_N$ on $N$
associated to  $\theta(v) \in \mathfrak{g}$.
\end{lemma}

As a consequence of Equation (\ref{equivariance}) and 
Lemma \ref{infinitesimalequivariance} we have

\begin{lemma}\label{inversion}
\hfill
\begin{enumerate}
\item $\Theta(hx) = h^{-1} \Theta(x), h \in  H, x \in M.$
\item $\Theta_{\ast}(V_M) = -V_N, v \in \mathfrak{h}.$
\end{enumerate}
\end{lemma}

We now prove 

\begin{lemma}\label{isometry}
$\Theta$ is a K\"ahler isomorphism.
\end{lemma}
\begin{proof}
Since $\theta:G \to G$ is holomorphic, any map of quotients it induces
is also holomorphic (by the universal property of quotients). 
Also any automorphism of a Lie algebra induces an isometry of the
Killing form (see \cite{Samelson}, pg. 14). Since in the semisimple
case the metric on $M$ is induced by the negative of the Killing form
on $\mathfrak{k}$ we are done in the semisimple case. For the case
of $GL_n(\C)$ we replace the Killing form by the trace form
and argue analogously.
\end{proof}

\subsection{The action of the Chevalley involution on the  momentum map}

The following result will tell us how $\Theta$ relates momentum levels
for symplectic quotients. Recall that  $N = M^{opp}$
so we have $\Theta:M \to N$.

\begin{proposition}\label{levels}
Let $\mu_M$ and $\mu_N$ be the momentum
maps for the actions of $T$ on $M$ and $N$ respectively.
Then there exists an element $\mathbf{\Lambda} \in (\mathfrak{t}^{\ast})^W$ such
that 
$$ \Theta^{\ast}\mu_N = \mathbf{\Lambda} - \mu_M.$$
\end{proposition}

\begin{proof}
For $v \in \mathfrak{t}$ we let $V_M$ and $V_N$ be the associated fundamental
vector fields on $M$ and $N$ respectively. Let $h^M_v$ and $h^N_v$ be the
Hamiltonian potentials of $V_M$ and $V_N$.  By Lemma \ref{inversion}
we have 
$$ \Theta_{\ast}(V_M) = - V_N.$$
We claim that there exists 
a linear functional $\mathbf{\Lambda} \in \mathfrak{t}^{\ast}$ such
$$\theta^{\ast} h^N_v = \mathbf{\Lambda}(v) - h^M_v.$$
To prove the claim it suffices to prove the differentiated version
$$\Theta^{\ast} dh^N_v = - dh^M_v.$$
Let $ p \in M$ and $w \in T_p(M)$. Then we have
\begin{align*}
\Theta^{\ast} dh^N_v|_p(w)  = 
\Theta^{\ast}(\iota_{V_N(\Theta(p))}\omega_N|_{\Theta(p)}(d\Theta|_p(w))) = 
\omega_N|_{\Theta(p)}(V_N(\Theta(p)),d\Theta|_p(w))\\ = 
-\omega_N|_{\Theta(p)}(d\Theta|_p(V_N(p)),d\Theta|_p(w))
 = -\omega_M(V_M(p),w) = - \iota_{V_M(p)}\omega_M(w).
 \end{align*}
The claim follows.

\medskip
It remains to prove that $\mathbf{\Lambda}$ is invariant under the Weyl group.
We first establish the $W$--equivariance of $\Theta^{\ast} \mu_N$.
To this end let $x \in M$ and $ w \in W$ and let $n(w) \in N(T)$ be
a representative of $w$ in $N(T)$, the normalizer of $W$ in $U(n)$.
We claim that we may choose these representatives such that $\theta(n(w)) = 
n(w)$. Indeed the Tits representatives, see \cite{MillsonToledano}, \S  2.5,
$\exp{(x_{\alpha})}\exp{(-x_{-\alpha})}
\exp{(x_{\alpha})}$ have this property because
$\exp{(x_{\alpha})}\exp{(-x_{-\alpha})}\exp{(x_{\alpha})} =
\exp{(-x_{-\alpha})}\exp{(x_{\alpha})}\exp{(-x_{-\alpha})}$ as can be checked
by a computation in $\mathfrak{sl}_2(\C)$.
We next  observe that it is an immediate consequence of the
$K$--equivariance of the momentum map for the action of $K$ on $M$
(and the relation between the $K$ and $T$ momentum maps) that
$$\mu_M(n(w)\cdot x) = Ad^{\ast}w(\mu_M(x)).$$
Then we have
$$\Theta^{\ast} \mu_N( n(w)\cdot x) = \mu_N(\Theta (n(w)\cdot x)) = 
\mu_N(n(w) \Theta(x))
= Ad^{\ast}w \mu_N(\Theta(x)) = Ad^{\ast}w (\Theta^{\ast} \mu_N)(x).$$ 
Since $\mu_M$ is also $W$--equivariant we find that
$\mathbf{\Lambda} = \Theta^{\ast} \mu_N + \mu_M$ is also $W$--equivariant.
Thus, since $\mathbf{\Lambda}(x)$ is a constant function 
 $\mathbf{\Lambda} = \mathbf{\Lambda}(wx) = Ad^{\ast}w(\mathbf{\Lambda}(x)) 
= Ad^{\ast}w (\mathbf{\Lambda})$. 

\end{proof}

\begin{corollary}
If $G$ is semisimple then $\mathbf{\Lambda} =0$ and we have
$$\Theta^{\ast} \mu_N = - \mu_M.$$
\end{corollary}

We obtain a general isomorphism formula for the action of $\Theta$
on symplectic quotients.
\begin{theorem}\label{generalsymplecticquotients}

The map $\Theta:M \to N$ induces a homeomorphism (K\"ahler
isomorphism in the smooth case)
$$\Theta: M//_{\br}T \to N//_{\mathbf{\Lambda} - \br}T.$$
\end{theorem}

\subsection{Formulas for the Chevalley involution}\label{formulas}

The above isomorphisms will be more useful if we have more 
explicit formulas for $\Theta$. We begin with a remarkably
useful lemma. It will apply to all simple complex groups except
$SL_n(\C)$ and $E_6$ (the case of $SO_{4n+2}(\C)$ will require a slight
modification, see below).

\begin{lemma}\label{inner}
Suppose the Chevalley involution is inner with $\theta = AdJ$. Let
$P$ be a parabolic subgroup and $M = G/P$. Then  
$\Theta:M \to M$ coincides with the map induced by action of $J$
by left translation.
\end{lemma}

\begin{proof}
The lemma is obvious if we use the description of $\Theta$ from 
Remark \ref{anotherdescription}. If $\theta = Ad J$ then we may take
$J$ as a representative of $n(w_0)$ and we obtain
$\Theta(gP) = AdJ(g) J P = Jg P.$
\end {proof}

Assume first that $G$ is a classical group other than $SL_n(\C)$ or
$SO_{4n+2}(\C)$. In what follows $J$ (or $J_G$) will denote an element of $G$
such that $\theta(g) = AdJ (g)$. 
Let $Q_n$ denote the $n$ by $n$ matrix with $1$'s on the counter-diagonal and 
$0$'s elsewhere.  
We have 
\[
J_{Sp_{2n}(\C)} = 
\begin{pmatrix}
0   & Q_n \\
- Q_n & 0
\end{pmatrix}, \quad
J_{SO_{4n}(\C)} = Q_{4n}, \quad 
J_{SO_{2n+1}(\C)} = (-1)^n \, Q_{2n+1}. 
\]
 
In the case of $SO_{4n+2}(\C)$ there is an element again denoted $J$ in 
$O_{4n+2}(\C)$ such that $\theta(g) = Ad J(g)$.  Here $J$ is simply $Q_{4n+2}$.

Note that in all cases $J$ is a scalar multiple of 
the matrix of the bilinear form relative to the standard basis. 

\subsection{Computations for  $GL_n(\C)$}

\hfill

We begin by relating the duality maps $\Psi$ and $\hat{\Psi}$
of the Introduction to the maps $\Theta$ and $\hat{\Theta}$
as promised in the Introduction.

\begin{lemma}\label{relationtoChevalleyinvolution}
\hfill
\begin{enumerate}
\item $\Theta = \Psi.$
\item $\hat{\Theta} = \hat{\Psi}.$
\end{enumerate}
\end{lemma}
\begin{proof}
We first give the proofs for Grassmannians.
To prove $(i)$ first note that the $j$-th column $C_j(\theta(g))$ is
the $j$--th row of $g^{-1}$ whence
$$(C_j(\theta(g)), C_i(g)) = \delta_{ij}.$$
Here $(\ ,\ )$ denotes the form $B$.
Hence the last $l= n-k$ columns of $\theta(g))$ are orthogonal
to the first $k$ columns of $g$ and we have a commutative diagram

\begin{center}
\begin{math}
\begin{CD}
G/P @>{\Theta}>>G/P^{opp} \\
@V{\pi_k}VV                  @VV{\pi_l}V \\
Gr_k(\C^n) @>{\Psi}>> Gr_l(\C^n) 
\end{CD}
\end{math}
\end{center} 

Here the vertical arrows are given by $\pi_k(gP) = g\cdot e_1 \wedge \cdots 
\wedge e_k$ and $\pi_l(gP^{opp}) = g\cdot e_{k+1}\wedge \cdots \wedge e_n$.
the first statement follows.

To prove the second statement we have only to observe there is another
commutative square.

\begin{center}
\begin{math}
\begin{CD}
G \times _P \mathbb{C} @>{\hat{\Theta}}>>G \times _{P^{opp}} \mathbb{C} \\
@V{f_k}VV                  @VV{f_l}V \\
\mathcal{T}_k @>{\hat{\Psi}}>> \mathcal{T}_l
\end{CD}
\end{math}
\end{center} 

Here the vertical arrows are given by $f_k([g,z]) = z g\cdot e_1 \wedge 
\cdots \wedge e_k$ and $f_l([g,w]) = w g \cdot e_{k+1}\wedge \cdots \wedge e_n$.
The bundle on the upper left is obtained from the equivalence relation
$(g,z) \sim (gp, det_k^{-1}(p) z)$ and the bundle on the upper right
is obtained from the equivalence relation $(g,w) \sim (gp, det_l^{-1}(p) w)$.

The reader will verify that the statements in the  lemma may now be deduced
for the case of flag manifolds for $GL_n(C)$ by comparing the diagram

\begin{center}
\begin{math}
\begin{CD}
F_\bk(\C^n) @>{\Psi}>> F_\bl(\C^n) \\
@V{i}VV                  @VV{i}V \\
\prod_{i \leq m} Gr_{k_i}(\C^n) 
@>F \circ {\prod_i \Psi_i}>> 
\prod_{i \leq m} Gr_{l_i}(\C^n) 
\end{CD}
\end{math}
\end{center} 

with the analogous diagram where $\Psi$ is replaced by $\Theta$ and
the factors $\Psi_i$ are replaced by the factors $\Theta_i$.

\end{proof}

In order to apply Proposition \ref{generalsymplecticquotients} we need
to compute $\mathbf{\Lambda}$. We have seen that $\mathbf{\Lambda}= 0$
in the semisimple case. We now compute $\mathbf{\Lambda}$ for
$GL_n(\C)$.  

\begin{proposition}\label{theconstantforthegenerallineargroupprop}
\hfill

Suppose that $\mu_M$ takes values in the orbit corresponding to
the $n$--tuple of eigenvalues (arranged in weakly decreasing order)
$\mathbf{\lambda} = (\lambda_1,\ldots,\lambda_n)$ and
$\mu_N$ takes values in the orbit corresponding to the $n$--tuple
of eigenvalues $(\nu_1,\ldots, \nu_n)$. Then 
$$\mathbf{\Lambda} = (\lambda_1 + \nu _n) \varpi_n.$$
\end{proposition}

We will apply the Proposition to the case where $\lambda_n$ and
$\mu_n$ are both zero. We record the resulting formula in the form
we will use it.

\begin{corollary}\label{theconstantforthegenerallineargroup}
Suppose that $\mu_M$ takes values in  the orbit of 
$\sum_{i=1}^{n-1} a_i \varpi_{i}$
and $\mu_N$ takes values in the orbit of $\sum_{i=1}^{n-1} b_i \varpi_{i}$.
Then 
$$\mathbf{\Lambda} = (\sum_{i=1}^{n-1} a_i) \varpi_n.$$
\end{corollary}
\begin{proof}
With the assumptions of the corollary we have $\lambda_1 = \sum_{i=1}^{n-1} a_i$
and $\nu_n =0$.
\end{proof}

We need a preliminary lemma. Let $\mathbf{e}$ be the 
standard coordinate flag in $\C^n$ $\mathbf{e}= 
(\C e_1, \C e_1 \wedge e_2, \ldots, \C e_1\wedge \cdots \wedge e_{n-1})$.

\begin{lemma}\label{valueate} 
Suppose $M$ is the manifold of full flags in $\C^n$. Then we have
\begin{enumerate}
\item $\Theta(\mathbf{e}) = 
(\C e_n, \C e_{n-1} \wedge e_n, \ldots, \C e_2\wedge \cdots \wedge e_{n}).$
\item Suppose $\mu_M$ is normalized so that it takes values in the
orbit with smallest eigenvalue equal to zero. Then
$$\mu_M(\mathbf{e}) = \sum_{i=1}^{n-1} a_i \varpi_i.$$
\end{enumerate} 
\end{lemma}

\begin{proof}
The first statement follows from the first statement in Lemma 
\ref{relationtoChevalleyinvolution}.
\smallskip

The second statement follows from an explicit computation of the moment map
for the Grassmannian $Gr_k(\C^n)$ with the symplectic structure given
by embedding it as the orbit of $\varpi_k$.
The moment
map is given in \cite{GGMS} as follows.  For an $n$ by $k$ matrix
$A$, define $\mu_i([A])$, for $1 \leq i \leq n$, as
$$\mu_i([A]) = \frac{\sum_{i \in J} |\det A(J)|^2}{\sum_J |\det
A(J)|^2},$$ where $J$ ranges over the $k$-element subsets of
$\{1,\ldots,n\}$, and $A(J)$ is the $k$ by $k$ submatrix of $A$
whose rows are the rows of $A$ indexed by $J$.  Then $\mu =
(\mu_1,\ldots,\mu_n)$ is the moment map for the torus.
We see then that the value of the moment map at the standard coordinate
$k$--plane is $\varpi_k$. If we change the symplectic structure
to be the one corresponding to the orbit of $a\varpi_k$ then the
value at $\mathbf{e}$ will be $a \varpi_k$. Since the moment map for the full 
flag
manifold under the diagonal torus action is the sum of the
moment maps for each factor the second statement follows.
\end{proof}

\begin{remark} The point here is to check that the value $\varpi_k$
is attained at the flag $\mathbf{e}$.
\end{remark}

The Proposition is now a  consequence of the following lemma

\begin{lemma}
\hfill

We have the following identity of $\mathfrak{t}^{\ast}$--valued functions on 
$M$:
$$\mu_M + \mu_N \circ \Theta = 
(\lambda_{max} \circ \mu_M + \lambda_{min} \circ \mu_N) \varpi_n.$$
\end{lemma}
\begin{proof}
For ease of notation we prove only the special case where the flag
manifold is the manifold of full flags.
Note first that the difference between the left-hand side and the right-hand
side is invariant under the action of $\mathbb{R}^2$ that translates
$\mu_M$ by $s \varpi_n$ and translates $\mu_N$ by $t \varpi_n$. Hence
it suffices to prove the formula in the case that the 
$\lambda_{min} \circ \mu_M = 0$ and $\lambda_{min} \circ \mu_N = 0$.
Hence we may assume that $\mu_M$ and $\mu_N$ are as in the corollary.
We now compute both sides on the standard coordinate flag 
$\mathbf{e}$. This will determine $\mathbf{\Lambda}$.
By the first statement in Lemma \ref{valueate} we have
$$\Theta(\mathbf{e}) = \Psi(\mathbf{e}) =
(\C e_n, \C e_{n-1} \wedge e_n, \ldots, \C e_2\wedge \cdots \wedge e_{n}).$$
Hence $\Theta(\mathbf{e}) = n(w_0)(\mathbf{e})$ where $n(w_0)
$ is a representative in $N(T)$ for the longest element $w_0$ in the Weyl group. 
Hence we have
$$\mu_N(\Theta(\mathbf{e})) =\mu_N(n(w_0)(\mathbf{e}))= Ad^{\ast}w_0 
(\mu_N)(\mathbf{e})
= Ad^{\ast}w_0 (\sum_{i=1}^{n-1} b_i \varpi_i).$$
Note that the last equality followed from the second statement in 
Lemma \ref{valueate}. Thus $\mu_N(\Theta(\mathbf{e}))$ has {\em first} 
coordinate 
equal to zero. On the other hand, by the second statement in Lemma 
\ref{valueate}
we see that $\mu_M(\mathbf{e})$ has first coordinate equal to
$\sum_{i=1}^{n-1} a_i$.
Since the
sum of these two vectors has all components equal we conclude that
all components of the sum are equal to $\sum_{i=j}^{n-1} a_i$
whence $\mathbf{\Lambda} = (\sum_{i=1}^{n-1} a_i) \varpi_n$.
But it is immediate that the eigenvalues of $\mu_M$ are the sums
$\sum_{i=j}^{n-1} a_i$. Hence the largest eigenvalue of $\mu_M$ is
$\sum_{i=1}^{n-1} a_i$ and since the smallest eigenvalue of 
$\mu_N$ is zero by definition the lemma is proved.
\end{proof}

\section{The Mumford Quotient}

\subsection{Definition of Mumford Quotient}\label{Mumford_Defn}

We refer the reader to \cite{Dolgachev} for
additional details.  Suppose that $G$ is a reductive algebraic group,
$V$ is a projective variety, and $\eta : G \times V \to V$ is 
regular action of $G$. Let
$\pi : \mathcal{L} \to V$ be an ample line bundle over $V$.
A $G$--linearization
of $\mathcal{L}$ is a regular action $\widetilde{\eta} : G \times
\mathcal{L} \to \mathcal{L}$ which is linear on fibers and 
makes the following diagram
commute:
\[
\begin{CD}
   G \times \mathcal{L} @>\widetilde{\eta}>> \mathcal{L} \\
   @V{id \times \pi}VV                 @VV{\pi}V \\
   G \times V @>\eta>> V
\end{CD}
\]

Given such a linearization, we automatically get linearizations on all
tensor powers $\mathcal{L}^{\otimes N}$ of $\mathcal{L}$.
Thus
$G$ has an action on sections $s$ of $\mathcal{L}^{\otimes N}$ given
by $(g \cdot s)(x) = g \cdot s(g^{-1} \cdot x) =
\widetilde{\eta}(g,s(\eta(g^{-1},x)))$.
Let $\Gamma(V,\mathcal{L}^{\otimes N})^G$ denote the $G$--invariant
holomorphic sections of $\mathcal{L}^{\otimes N}$.
We define the semistable points of $V$ for the chosen linearization
$\widetilde{\eta}$ to be
$$ V_{\widetilde{\eta}}^{ss} = \{x \in V \, | \, (\exists N)
(\exists s \in \Gamma(V,\mathcal{L}^{\otimes N})^G)(s(x) \neq 0)\}.$$
The Mumford quotient $V //_{\widetilde{\eta}}\: G$
is defined as the quotient space of $V_{\widetilde{\eta}}^{ss}$ such that two
points $x,y \in V_{\widetilde{\eta}}^{ss}$ are identified iff their
$G$--orbit closures (computed in $V_{\widetilde{\eta}}^{ss}$) 
$cl(G \cdot x)$ and $cl(G \cdot y)$ intersect
nontrivially.  The Mumford quotient $V //_{\widetilde{\eta}} \, G$
is then a projective variety corresponding to the geometric points of
$\t{Proj}(S^G)$ where $S^G$ is the graded ring $\bigoplus_{N \geq
0} \Gamma(V,\mathcal{L}^{\otimes N})^G$.

In the case where we have an action $\eta$ of the complex torus 
$H \subset G = GL_n(\C)$, 
and a homogeneous line bundle $\mathcal{L} = G \times_P \C$ over $G/P$,
the set
of possible linearizations correspond to the complex characters of
$H$, which are all of the form $(z_1,...,z_n) \mapsto
\prod_{i=1}^n z_i^{r_i}$ where $\br = (r_1,...,r_n) \in \Z^n$.
(See \cite{Dolgachev}, Chapter 7.)
We denote the linearization associated to $\br$ as
$V //_\br \, H$.

\subsection{Gel'fand-MacPherson duality of Mumford quotients}

In the previous section the Mumford (or categorical) quotient
was given for a group acting on a projective variety, but there is an 
easier definition when 
the variety $V$ is affine.  We may take  
the trivial line bundle $\mathcal{L}$ where the sections
are simply $\mathcal{O}(V)$.  Any character $\chi$ of $G$ determines a
linearization of the trivial line bundle.  
Semistability is defined as
before.

Let $k = m+1$ and let 
$\mathcal{L}$ be the trivial line bundle $\C^{n \times k}
\times \C \rightarrow \C^{n \times k}$.  The group $GL_k(\C)$
acts on the right of $\C^{n \times k}$ by matrix multiplication.  
The group $H$ of nonsingular diagonal $n$ by $n$ complex matrices 
acts on the left of $\C^{n \times k}$.

Let the character $det^a : GL_k(\C) \rightarrow \C^\ast$ and the
character $\chi_\br : H \rightarrow \C^\ast$ be given.  The 
one--dimensional subgroup $K = \{(zI_n,z^{-1}I_k) : z \in \C^\ast\}$ of 
$H \times GL_k(\C)$ acts trivially on $\C^{n \times k}$.  
Let $G$ be the quotient of $H \times GL_k(\C)$ by $K$.  
The character $\chi_\br \times det^a$ defines a character of $G$ iff 
$|\br|=ak$, and we assume that is the case so that we have a
$G$--linearization of the trivial line bundle.  

The quotient by $GL_k(\C)$ alone is the 
Grassmannian $Gr_k(\C^n)$.  
There is a canonical line bundle $\mathcal{L}_1$ over $Gr_k(\C^n)$
which is the quotient of the pullback of $\mathcal{L}$ over the
semistable points.  The character $\chi_\br$ now
defines an $H$--linearization of $\mathcal{L}_1$.   
If we now take a Mumford quotient of
$Gr_k(\C^n)$ by $H$ using the character $\br$ this results in the
quotient $\C^{n \times k} //_{\br \times det^a} H \times GL_k(\C)$.

On the other hand, if each $r_i > 0$ then the 
quotient by $H$ alone is $(\C P^{k-1})^n$.  
Again there is a canonical quotient bundle $\mathcal{L}_2$ of the
trivial bundle pulled back over the semistable points, 
and the character $det^a$ determines a
$GL_k(\C)$--linearization of $\mathcal{L}_2$.  
If we now take the Mumford
quotient by $GL_k(\C)$ (projective equivalence) using the character
$det^a$ we again must get 
$\C^{n \times k} //_{\br \times det^a} H \times GL_k(\C)$.

This is illustrated by the following diagram:

\[
\begin{CD}
   \C^{n \times (m+1)}  @>{//_{det^a} GL_{m+1}(\C)}>>  Gr_{m+1}(\C^n) \\
   @V{//_\br H}VV   @VV{//_\br H}V \\
   (\C P^m)^n @>{//_{det^a} GL_{m+1}(\C)}>> \Mu
\end{CD}
\]
where $\Mu$ denotes the Mumford quotient $Gr_k(\C^n)//_\br H$. 
 
\section{Duality for torus quotients of Grassmannians on the quantum level}

\subsection{The Hodge star operator}

To promote the duality map $\Psi$ to a map of bundles we
recall the definition of the complex Hodge star operator $\ast$. Choose an
orientation for $\mathbb{C}^n$ (whence an orientation for 
$(\mathbb{C}^n)^{\ast})$. 
Let  $vol \in \bigwedge^n (\mathbb{C}^n)^{\ast}$ be the positively oriented 
element of
unit length for the form induced by $B$. We will say $vol$ is a complex
volume form. The complex volume form induces a map (by interior multiplication) 
$\alpha:\bigwedge^k(\mathbb{C}^n) \to \bigwedge^{n-k} (\mathbb{C}^n)^{\ast}$.
The form $B$ induces a map $\beta:\bigwedge^{n-k} (\mathbb{C}^n)^{\ast} \to
\bigwedge^{n-k} (\mathbb{C}^n)$. We define the complex Hodge star $\ast$
to be the composition $\beta \circ \alpha$. We note that $\ast$ is
an invertible linear map from $\bigwedge^k(\mathbb{C}^n)$ to
$\bigwedge^{n-k}(\mathbb{C}^n)$. The pair $\Psi$ and $\ast$ induce a map from 
the trivial $\bigwedge^k(\mathbb{C}^n)$--bundle over $Gr_{k}(\mathbb{C}^n)$ 
to the trivial $\bigwedge^{n-k}(\mathbb{C}^n)$--bundle over 
$Gr_{n-k}(\mathbb{C}^n)$. This bundle map carries the subbundle  
$\mathcal{T}_k$ to the subbundle $\mathcal{T}_{n-k}$. 
Consequently it induces a bundle isomorphism
from the tautological line bundle 
$\mathcal{T}_k$ over $Gr_{k}(\mathbb{C}^n)$
to  the tautological line  bundle $\mathcal{T}_{n-k}$ over 
$Gr_{n-k}(\mathbb{C}^n)$
covering $\Psi$ (so technically a bundle isomorphism from $\mathcal{T}_k$ to
$\Psi^{\ast} \mathcal{T}_{n-k})$. We obtain the following diagram

\begin{center}
\begin{math}
\begin{CD}
\mathcal{T}_k @>>> Gr_{k}(\mathbb{C}^n) \times \bigwedge^k(\mathbb{C}^n)\\
@V{\Psi \times\ast}VV                  @VV{\Psi \times\ast}V \\
\mathcal{T}_{n-k} @>>> Gr_{n-k}(\mathbb{C}^n) \times 
\bigwedge^{n-k}(\mathbb{C}^n)
\end{CD}
\end{math}
\end{center}

In order to obtain an isomorphism
of their duals we dualize the definition of $\ast$ to obtain
a new isomorphism again denoted $\ast$ from $\bigwedge^k((\mathbb{C}^n)^{\ast})$ 
to
$\bigwedge^{n-k}((\mathbb{C}^n)^{\ast})$. We obtain an induced isomorphism
of line bundles homomorphisms $\hat{\Psi}$ 
covering $\Psi$ from $\mathcal{L}_k$ to $\mathcal{L}_{n-k}$ by
dualizing the above diagram. 

\subsection{The Homogeneous Coordinate Ring of $Gr_k(\C^n) //_\br \, H$}

In this section we describe a basis for the coordinate ring of $Gr_k(\C^n)
//_\br \, H$. We begin by recalling that the Pl\"ucker
embedding $\iota_k$ of $Gr_k(\C^n)$ is the projective embedding of
$Gr_k(\C^n)$ corresponding to the very ample line bundle $\mathcal{L}$.
According to the general theory of projective embeddings and
line bundles we have an embedding
$$\iota_k:Gr_k(\C^n) \to \mathbb{P}(\Gamma(Gr_k(\C^n),\mathcal{L})^{\ast}).$$
It is standard that $\Gamma(Gr_k(\C^n),\mathcal{L}) \cong
\bigwedge^k((\C^n)^{\ast})$. In what follows we will need an explicit formula
for this isomorphism. 

\vskip 12pt

Let $x \in Gr_k(\C^n)$ and $\tau \in \bigwedge^k((\C^n)^{\ast})$. We let
$res_x:\bigwedge^k((\C^n)^{\ast}) \to \bigwedge^k((x)^{\ast})= \mathcal{L}_x$
be the operation of restriction of covectors to $x$. If 
$\tau \in \bigwedge^k((\C^n)^{\ast})$ we let $\widetilde{\tau}$ be
the section of $\mathcal{L}$ defined by
$$\widetilde{\tau}(x) = res_x(\tau).$$
The following lemma is then standard

\begin{lemma}\label{globalsections}
The map $\tau \rightarrow \widetilde{\tau}$ induces an isomorphism
$\bigwedge^k((\C^n)^{\ast}) \cong \Gamma(Gr_k(\C^n),\mathcal{L}).$
\end{lemma}
Let $\theta_i, 1 \leq i \leq n$ be the basis for $(\C^n)^{\ast}$
dual to the standard basis $\epsilon_i, 1 \leq i \leq n$. For
$I = \{i_1,..,i_k\} \subset \{1,2,...,n\}$ with $i_1 < i_2 < \cdots i_k$
we define
$$\theta_I = \theta_1 \wedge \theta_2 \wedge \cdots \wedge \theta_k.$$
The linear functions $\theta_I$ as $I$ ranges through the $k$--element 
subsets of $\{1,2,...,n\}$ give a basis for $\bigwedge^k((\C^n)^{\ast})$ and 
consequently
give homogeneous coordinates (the Pl\"ucker coordinates) to be denoted
$X_{i_1,i_2,...,i_k}$ (or $X_I$) on the projective space
$\mathbb{P}(\bigwedge(\C^n))$,

\subsubsection{A basis for the homogeneous coordinate ring of $Gr_k(\C^n)$}
\label{rectstandardbasis}
We begin by noting that we have the following formula for
the homogeneous coordinate ring $R_k$ of $Gr_k(\C^n)$ as a graded vector
space.

\begin{equation*}\label{gradedring}
R_k = \bigoplus_{N=0}^{\infty} \Gamma(Gr_k(\C^n), \mathcal{L}^{\otimes N})
= \bigoplus_{N=0}^{\infty} V_{N\varpi_k}.
\end{equation*}

Here $ V_{N\varpi_k}$ is the irreducible representation of $GL_n(\C)$
with highest weight $N\varpi_k$. 
Let $\widetilde{R}_k = \bigoplus_{N=0}^{\infty} R_k^{(N)}$ be the graded ring 
$\b{C}[X_{i_1,...,i_k}]$,
where $i_1,...,i_k$ ranges over $k$-element subsets of
$\{1,...,n\}$. We have seen that $\bigwedge((\C^n)^{\ast}) \cong
\Gamma(Gr_k(\C^n), \mathcal{L}))$ whence the Pl\"ucker coordinates
correspond to a basis of the degree one elements of $R_k$.
We obtain a map of rings $\phi:\b{C}[X_{i_1,...,i_k}] \to R_k$.
The following lemma is standard. 

\begin{lemma}

 \hfill

\begin{enumerate}
\item $\phi$ is onto.
\item The kernel of $\phi$ is generated by quadratic relations in the
Pl\"ucker coordinates called the Pl\"ucker relations.
\end{enumerate} 
\end{lemma}

\vskip 12pt
 
We now 
define certain elements $f_T \in R_k$ defined by fillings $T$ of
rectangular Young diagrams $D$ with $k$ rows by numbers between $0$ and $n$. 
Let $D_N$ be the rectangular Young diagram with $k$ rows and $N$ columns and 
let $T$ be a filling of $D_N$.
Let $I_i$ be the entries in the $i$--th column of $T$. Define
$$f_I = X_{I_1}X_{I_2}\cdots X_{I_N}.$$
We define $deg(T)$, the degree of $T$,
by 
$$deg(T) = N$$
and note that $f_T $ is in the $N$-th graded summand $R_k^{(N)}$ of $R_k$. 
Thus we find 
$deg(f_T)$, the degree of $f_T$ relative to the above grading, is also $N$.
We can now describe a basis for $R_k^{(N)}$.
We recall that a filling $T$ of $D_N$ is said to be semistandard if
the entries in each column are strictly increasing and the entries in 
each row are weakly increasing. 
We will use $\mathcal{SS}(D,n)$ to denote the set of semistandard
fillings of a Young diagram $D$ by the integers $1,2,...,n$.
We have, see \cite{DO},Chapter I, 
Theorem 1.

\begin{theorem} \label{standardbasis}
The functions $f_T$ as $T \in \mathcal{SS}(D_N,n)$ form a basis for $R_k^{(N)}$.
\end{theorem}

We will call the basis of $R_k$ given by the set of $f_T$'s with
$T$ semistandard, the standard basis of $R_k$.

\vskip 12pt

For $i$ between $1$ and $n$ we define $w_i(T)$ to be
the number of times $i$ appears in $T$ and we define the weight $wt(T)$ of $T$
to be the $n$-tuple $wt(T) = (w_1, \cdots w_n)$. 
We also define the weight $wt(f_T)$ by $wt(f_T) = wt(T)$. This terminology
is justified by the following

\begin{lemma}\label{weightbasis}
Under the action of $H$ on the graded ring $R_k$ the function $f_T$
is a weight vector of weight $wt(T)$.
\end{lemma}

\subsubsection{A basis for the homogeneous coordinate ring $S_k$ of 
$Gr_k(\C^n)//_{\br}H$}
Let $\br =
(r_1,...,r_n)$ be a tuple of non-negative integers.  Let
$a = \frac{\sum_{i=1}^n r_i}{k}$ as usual. We assume that
$a$ is an integer. We have
the following formula analogous to Equation \ref{gradedring} for
the underlying graded vector space of $S_k$

\begin{equation*}\label{subgradedring}
S_k = \bigoplus_{N=0}^{\infty} 
\Gamma(Gr_k(\C^n), \mathcal{L}^{\otimes Na}(N\br))^H
= \bigoplus_{N=0}^{\infty} V_{N\varpi_k}(N\br).
\end{equation*}

Note that if $f_T \in S_k^ N$ then we have

\begin{equation*}\label{weightanddegree}
deg(T) = N a \ \text{and} \ wt(T) = N \br.
\end{equation*}

We now check that this relation between $\br$ and the degree and weight of $T$
is automatically satisfied if we merely assume that $wt(T)$ is
an integral multiple of $\br$. We thereby obtain a simpler description
of the subring of $H$--invariants $S_k \subset R_k$.

\begin{lemma}
$S_k$ is the subring of $R_k$
spanned by the monomials $f_T$ such that
$wt(T)$ is a multiple of $\br$.
\end{lemma}

\begin{proof}
Suppose that $wt(T) = \ell \br$ with $\ell$ a positive integer.
Suppose $deg(T) = M$. Then $T$ is a filling of the $k$ by $M$
rectangle $D$ whence (equating the total number of boxes in $D$)
$$k M = wt(T) = \ell |\br| \ \text{so} \ deg(T) = M = \ell a.$$
\end{proof}

As a consequence of the previous lemma, Theorem \ref{standardbasis}
and Lemma \ref{weightbasis} we obtain a basis for $S_k$.

\begin{theorem}
The set of standard basis vectors $f_T$ with weight a multiple
of $\br$ is a basis for $S_k$.
\end{theorem}
We will call the resulting basis the standard basis of $S_k$.

\subsection{The proof of Theorem \ref{rings}}\label{ringtheorem}
In this subsection we prove Theorem \ref{rings}.
We begin by proving the $H$--equivariance of $\ast$ and the bundle
map $\hat{\Psi}$. First we deal with $\ast$.

\begin{lemma}
Let $u \in \bigwedge^k(\mathbb{C}^n)$ and $g \in GL_n(\C)$.
Then
$$\ast g \ u = det(g) (g^t)^{-1} \ast u.$$
In particular we have
$$\ast h \ u = det(h) h^{-1} \ast u, h \in H.$$
\end{lemma}
\begin{proof}
We recall from the Introduction that $\ast = \beta \circ \alpha$
where $\alpha$ is given by contraction with the volume form $vol$
and $\beta:\bigwedge^{n-k}((\mathbb{C}^n)^{\ast}) \to 
\bigwedge^{n-k}((\mathbb{C}^n)$ is the map induced by $B$.
We prove suitable equivariance formulae for each of $\alpha$ and $\beta$.
First we claim that for any $g \in GL_n(\mathbb{C})$ we have
$$\alpha(gu) = det(g) (g\alpha)(u).$$
Indeed let $v \in \bigwedge^{n-k}((\mathbb{C}^n)$
\begin{align*}
\alpha(gu) = vol(gu,v) = vol(gu,gg^{-1}v) = det(g) vol(u, g^{-1}v)
=det(g)(g\alpha)(v)
\end{align*}
We conclude by observing that for $g \in GL_n(\C)$ and $\tau \in 
\bigwedge^{n-k}((\mathbb{C}^n)^{\ast})$
we have
$$\beta(g \tau) = (g^t)^{-1}\beta(\tau).$$
Indeed it suffices to prove that for $v \in \mathbb{C}^n$ we
have $b(gv) = (g^t)^{-1}b(v)$ where $b:\mathbb{C}^n \to (\mathbb{C}^n)^{\ast}$
is the map induced by the bilinear form $B$. But this is immediate.
\end{proof}

\vskip 12pt

The next corollary follows by dualizing the previous lemma.

\begin{corollary}
Let $\eta \in \bigwedge^k((\mathbb{C}^n)^{\ast})$ and $h \in H$.
Then we have
$$\ast g \ \eta = det(g)^{-1} (g^t)^{-1} \ast \eta.$$
In particular we have
$$\ast h \ \eta = det(h)^{-1} h^{-1} \eta.$$
\end{corollary}

We need another corollary. Let $\rho_k$ be the $k$--the exterior power of
the standard representation of $GL_n$ and define $\rho_k^{\theta}$
by $\rho_k^{\theta} = \rho_k \circ \theta$. Recall that $\theta$
is the Chevalley involution $\theta(g) = (g^t)^{-1}$.

Then we have by the above (since $\widetilde{\Psi_k}$ is equal to the Hodge star
on $\bigwedge^k((\C^n)^{\ast}))$

\begin{corollary}\label{intertwining1}

$$\widetilde{\Psi}_k\circ \rho_k \circ \widetilde{\Psi}_{n-k} = \rho_k^{\theta} \otimes
det^{-1}.$$

\end{corollary}

Next we deal with $\hat{\Psi}:\mathcal{L}_k^{\otimes m a} 
\to \mathcal{L}_{n-k}^{\otimes m a}$ and prove the first statement in
Theorem \ref{rings}.

\begin{lemma}\label{equivarianceofthebundlemap}
$$\hat{\Psi} \circ h = h^{-1} \circ \hat{\Psi}.$$
\end{lemma}

\begin{proof}
It suffices to prove the lemma for the case $m=1$.
We have 
$$h(x,\alpha^{\otimes a}) = 
(hx, \chi_{\br}(h) (h\alpha)^{\otimes a})$$
whence
\begin{align*}
\hat{\Psi}(h(x,\alpha^{\otimes a})) =  (\Psi(hx), \chi_{\br}(h)
(\ast h \alpha)^{\otimes a})  = & \ 
(h^{-1}\Psi(x), \chi_{\br}(h)(det(h)^{-1}
h^{-1}\ast \alpha)^{\otimes a}) \\
=(h^{-1}\Psi(x), \chi_{\br}(h) det(h)^{-a} h^{-1}(\ast \alpha)^
{\otimes a}) = & (h^{-1}\Psi(x), \chi_{\L -\br}(h^{-1})
h^{-1}(\ast \alpha)^ {\otimes a}) =  h^{-1}
 \hat{\Psi}((x,\alpha^{\otimes a}).
 \end{align*}
 \end{proof}

Now we prove the second statement in Theorem \ref{rings}. 
Note that the map $\widetilde{\Psi}$ on sections induced by the bundle map 
$\hat{\Psi}$ is given by
$$\widetilde{\Psi}(s)(x) = \hat{\Psi}(s(\Psi^{-1}(x))).$$
The second statement in Theorem \ref{rings} follows from
\begin{lemma}
A section $s$ is $H$--invariant $\Leftrightarrow$ the section
$\widetilde{\Psi}(s)$ is $H$--invariant.
\end{lemma}

\begin{proof}
By symmetry it suffices to prove the direction $\Rightarrow$.
So assume that $s$ is $H$--invariant. It suffices to prove
that $h^{-1}\widetilde{\Psi}=\widetilde{\Psi}$ for all $h\in H$.
But
\begin{align*}
&(h^{-1}\widetilde{\Psi})(x) =  h^{-1}\widetilde{\Psi}(hx) = 
h^{-1}\hat{\Psi}(s(\Psi^{-1}(hx))) =h^{-1} \hat{\Psi}(s(h^{-1}(\Psi^{-1}(x))))\\
& =\hat{\Psi}(h(s(h^{-1}(\Psi^{-1}(x))))) =  \hat{\Psi}(s(\Psi^{-1}(x)))
=\widetilde{\Psi}(s)(x).
\end{align*}
\end{proof}

Later we will need the following immediate consequence of 
Corollary \ref{intertwining1}.
\begin{lemma}\label{intertwining2}
Let $\rho_{N,a,k}$ be the representation of $GL_n(\C)$
on the vector space of sections 
$\Gamma(Gr_k(\C^n), \mathcal{L}_k^{\otimes aN})$.
Let $\rho_{N,a,k}^{\theta} = \rho_{N,a,k} \circ \theta$. 
Then we have
$$\widetilde{\Psi}_{N,a,k} \circ \rho_{N,a,k} \circ \widetilde{\Psi}_{N,a,n-k} = 
\rho_{N,a,n-k}^{\theta}.$$
\end{lemma}

\subsection{The proof of Theorem \ref{tableaux}}
In this subsection we prove the formula of the introduction for the
action of the isomorphism $\widetilde{\Psi}$ on the standard basis
vectors $f_T$ for the homogeneous coordinate ring of $Gr_k(\C^n)$.

\vskip 12pt

The bundle map  $\hat{\Psi}$ satisfies  the formula
$$ \hat{\Psi}(res_x(\tau)) = res_{\Psi(x)}(\ast(\tau)), \quad
\tau \in \bigwedge ^k((\mathbb{C}^n)^{\ast})$$
We now have
\begin{lemma}
$$\hat{\Psi}(\widetilde{\theta_I}) = sgn(I,J) \ \widetilde{\theta_J}.$$
\end{lemma}
\begin{proof} Let $x \in Gr_k(\C^n)$.
Then we have $\widetilde{\theta_I}(x) = res_x({\theta_I})$.
Now let $y \in Gr_{n-k}(\C^n)$. We have 
\begin{align*}
\widetilde{\Psi}(\widetilde{\theta_I}(y))= 
\hat{\Psi}(\widetilde{\theta_I}(\Psi^{-1}y)) = 
res_{\Psi(\Psi^{-1}(y))} (\ast \theta_I) = sgn(I,J) \ res_y (\theta_J)
=sgn(I,J) \  \widetilde{\theta}(y).
\end{align*}
\end{proof}

The following corollary is an immediate consequence of the lemma.
 
 \begin{corollary}
 $$\widetilde{\Psi}(\widetilde{\theta_{I_1}}\otimes \widetilde{\theta_{I_2}}
 \otimes \cdots \otimes \widetilde{\theta_{I_p}})= sgn(I_1,J_1)\cdots 
 sgn(I_p,J_p) \ 
 \widetilde{\theta_{J_1}}\otimes \widetilde{\theta_{J_2}}
 \otimes \cdots \otimes \widetilde{\theta_{J_p}}.$$
 \end{corollary}
 
Theorem \ref{tableaux} is a consequence of the corollary because $f_T$ is the
image of a tensor product of the above form under the linear map given
by multiplication of sections.

\section{Duality for torus quotients of  flag manifolds on the quantum level}
\label{Flag}
\subsection{The relation between partial flag manifolds and products of 
Grassmannians}
In this section we extend our results on duality of homogeneous coordinate rings
of torus quotients of Grassmannians to torus quotients of flag manifolds.
\subsection{The homogeneous coordinate ring of the flag manifold}
Let $P$ be a parabolic subgroup of $G = \t{GL}_n(\C)$.  
Then $G/P$ is a partial flag manifold $F_{k_1,...,k_m}(\C^n)$ where $0 < k_1 < 
k_2 <
\cdots < k_m < n$.  The ample line bundles over $G/P$ are parametrized
by weights $a_1 \varpi_{k_1} + \cdots + a_m \varpi_{k_m}$ where each
$a_i$ is a positive integer.  
In what follows, let $\bk = (k_1,...,k_m)$ and $\ba = (a_1,...,a_m)$
and we will abbreviate the above dominant weight to $\lambda_{\ba}$.
Choose an $m$--tuple $\ba$ as above. The $m$-tuple $\ba$ corresponds to
a character $\chi_{\ba}$ of $P$. We let $\Lf$ 
be the line bundle over $\Fk$ with isotropy representation $\chi_{\ba}^{-1}$
(so $\Lf$  had total space defined by the equivalence relation 
$(gp,\chi_{\ba}(p) z) \sim (g,z)$). We define 
$V_{\lambda_{\ba}}=\Gamma(\Fk, \mathcal{L}_{\ba})^{\ast}$. The group $G$
acts on $V_{\lambda_{\ba}}$ and the  flag manifold $F_{k_1,...,k_m}(\C^n)$
is embedded in $\mathbb{P}(V_{\lambda_{\ba}})$ as the orbit of the line
through a highest weight vector.

In this section we will use the embedding 

$$i: \Fk \to Gr_{k_1}(\C^n) \times \cdots \times Gr_{k_m}(\C^n)$$

to promote our duality results for Grassmannians to flag manifolds.

The line bundle $\Lf$ is very ample and we
obtain an equivariant projective embedding
$$\iota:\Fk  \to \mathbb{P}(\Gamma(\Fk, \Lf)^{\ast}) =
\mathbb{P}((V_{\lambda_{\ba}})^{\ast}).$$

Accordingly we have the  following formula for the homogeneous coordinate ring
$\Rk$ of $\Fk$
$$\Rk = \bigoplus_{N=0}^{\infty}\Gamma(\Fk, (\Lf)^{\otimes N})
= \bigoplus_{N=0}^{\infty} V_{N\lambda_{\ba}}.$$

\vskip 12pt

We have the very ample line bundle
$\mathcal{L}_{k_1}^{a_1} \boxtimes \cdots \boxtimes
\mathcal{L}_{k_m}^{a_m}$ over the product $Gr_{k_1}(\C^n)
\times \cdots \times Gr_{k_p}(\C^n)$.

We will also use $\widetilde{\Rk}$ to denote the homogeneous coordinate
ring of $Gr_{k_1}(\C^n)\times \cdots \times Gr_{k_p}(\C^n)$. Hence
$$\widetilde{\Rk} = \bigoplus_{N=0}^{\infty}
\Gamma( Gr_{k_1}(\C^n)\times \cdots \times Gr_{k_p}(\C^n),  
\mathcal{L}_{k_1}^{\otimes Na_1} \boxtimes \cdots \boxtimes
\mathcal{L}_{k_m}^{\otimes Na_m})
=V_{a_1 \varpi_1}\otimes \cdots \otimes 
V_{a_m \varpi_m}.$$

Recall that the irreducible representation 
$V_{\lambda_{\ba}}$ occurs with multiplicity one in the
tensor product $V_{a_1 \varpi_1}\otimes \cdots \otimes 
V_{a_m \varpi_m}$. Hence there is a canonical $GL_n(\C)$--quotient
mapping $\pi:V_{a_1 \varpi_1}\otimes \cdots \otimes 
V_{a_m \varpi_m} \to V_{\lambda_{\ba}}$. We let $\alpha =\pi^{\ast}$
be the dual map.

The following lemma will be very important in what follows. 

\begin{lemma}\label{pullback}
$$i^{\ast}(\mathcal{L}_{k_1}^{a_1} \boxtimes \cdots \boxtimes
\mathcal{L}_{k_m}^{a_m}) = \Lf.$$
\end{lemma}

\begin{proof}
We consider the following diagram

\begin{center}
\begin{math}
\begin{CD} 
F_\bk(\C^n) @>{\Delta}>> F_\bk(\C^n) \times \cdots  \times F_\bk(\C^n)\\
@V{Id}VV                  @VV{\pi}V \\
F_\bk(\C^n)   
@>{i}>> 
Gr_{k_1}(\C^n) \times \cdots \times Gr_{k_m}(\C^n)
\end{CD}
\end{math}
\end{center}

We need the following simple general observation whose proof we leave to the 
reader. 
Let $P \subset Q$ be subgroups of 
a group $G$. Consequently we have a projection $\pi: G/P \to G/Q$.
Let $\chi$ be a character of $Q$ and let $\mathcal{L}$ be the homogeneous 
line bundle with isotropy representation $\chi$, that is, the line bundle 
with total space $G \times _Q \C$ where $(g,z) \sim (gq, \chi(q)^{-1}z)$.
Then the pull-back of $\mathcal{L}$ to $G/P$ is the homogeneous line bundle
with isotropy representation $\chi|P$.

From this observation we find that 
$\pi^{\ast}(\mathcal{L}_{k_1}^{a_1} \boxtimes \cdots \boxtimes
\mathcal{L}_{k_m}^{a_m})$ is an outer tensor product of the same
form except the isotropy representations are the restrictions of the characters
corresponding to the weights 
$a_i\varpi_{k_i}, 1 \leq i \leq m$ to $P$. The pull-back of this outer tensor 
product by the
diagonal map $\Delta$ gives the inner tensor product. But the inner
tensor product of homogeneous line bundles is again homogeneous with
isotropy character the product of the characters of the factors.
But this product is just the character corresponding to the weight
$a_1 \varpi_{k_1} + \cdots + a_m \varpi_{k_m}$.
  
\end{proof}

Combining these observations we obtain the following
commutative diagram. 

\begin{center}
\begin{math}
\begin{CD}
F_\bk(\C^n) @>{\iota}>> \mathbb{P}((V_{\lambda_{\ba}})^{\ast}) \\
@V{i}VV                  @VV{\alpha}V \\
Gr_{k_1}(\C^n) \times \cdots \times Gr_{k_m}(\C^n)   
@>>> 
\mathbb{P}((V_{a_1 \varpi_1}\otimes \cdots \otimes 
V_{a_m \varpi_m})^{\ast})
\end{CD}
\end{math}
\end{center}

\vskip 12pt

\subsubsection{The $\br$-linearization of the action of $PH$}

Let $H$ be the complex torus of diagonal matrices in $GL_n(\C)$,
and let $PH$ be the image of $H$ under the quotient map 
$GL_n(\C) \to PGL_n(\C)$.
The $H$--linearizations of $\Lf$
correspond to a character $\chi_\br$ of $H$ given by 
$(z_1,...,z_n) \to \prod_i z_i^{r_i}$, where $\br \in \Z^n$.  
We denote the linearized bundle associated to this character as
$\Lf(\br)$. 

\begin{lemma}\label{existenceoflinearization}
The induced action of $H$ on $\Lf$ corresponding to
$\br$ descends to the 
quotient group $PH$ iff \  $|\br| = \sum_i a_i k_i$.   
\end{lemma}

\begin{proof}
Let $h = \mu I$ be a nonzero scalar matrix. Then for $[g,z]$ in the
total space of $\Lf(\br)$ we have
$$h [g,z] = [hg, \chi_{\br}(h)z] = [gh, \mu^{|\br|} z] = 
[g, (\prod_i det_{k_i}^{a_i}(\mu I_n))^{-1} \mu^{|\br|} z] = 
[g, \mu^{-\sum_i a_i k_i} \mu^{|\br|} z].$$ 
Thus $ h[g,z] = [g,z] \Leftrightarrow |\br| = \sum_i a_i k_i$.
\end{proof}

\subsubsection{The standard basis of $\Rk$}
In  what follows we will construct the standard basis of 
$\Gamma(\Fk, \Lf)$.  Since the basis vectors are
weight vectors for $H$ we will also obtain a basis
for $\Gamma(\Fk, \Lf)(\br)$. By replacing $\lambda_{\ba}$
by $\lambda_{N\ba}$ one obtains the standard basis for the $N$-th graded
summand of the homogeneous coordinate ring $\Rk$.

Define a partition $\mathbf{p}=(p_1,...,p_n)$ by
$$p_i = \sum_{j=i}^n a_j.$$
Let $D$ be the Young diagram corresponding to the partition $\mathbf{p}$
(so there are $p_i$ boxes in the $i$--th row, $1 \leq i \leq m$).
Let $T$ be a filling (not necessarily semistandard but such that the entries 
in each column are strictly
decreasing) of $D$ by elements in $1,2,...,n$. 
Our goal is to construct an element $f_T \in \Rk.$ 
The key to doing this is first to construct a section of the
line bundle
$\mathcal{L}_{k_1}^{a_1} \boxtimes \cdots \boxtimes
\mathcal{L}_{k_m}^{a_m}$ over the product $Gr_{k_1}(\C^n)
\times \cdots \times Gr_{k_p}(\C^n)$ then pull-back to $\Fk$
using Lemma \ref{pullback}.

To do this we divide $T$ up into $m$ rectangular
subtableaux $T_1,T_2,...,T_m$ where $T_{n+1 -i}$ is a  filling of 
a $k_i$ by $a_i$ rectangle. Thus we  take $T_1$ to be the rectangle
that is the union of the last $a_1$ columns. See Example \ref{subdivide}.

Next observe that $f_{T_i}$ is a section of the line bundle 
$\mathcal{L}_{k_i}^{a_i}, 1 \leq i \leq m$. 
The tensor product $f_{T_1} \otimes \cdots \otimes f_{T_m}$
is the desired section of $\mathcal{L}_{k_1}^{a_1} \boxtimes \cdots \boxtimes
\mathcal{L}_{k_m}^{a_m}$

We define $f_T$ by

$$f_T = i^{\ast}(f_{T_1} \otimes \cdots \otimes f_{T_m}).$$

Here $i^{\ast}$ denotes the pullback operation from sections
of a line bundle to sections of the pull-back bundle.

By Lemma \ref{pullback} we have
\begin{lemma}
$f_T$ is a section of $\mathcal{L}_{\bk}^{\ba}$.
\end{lemma}

We now have the following theorem, see \cite{GonciuleaLakshmibai},
Chapter 7, Theorem 2.1.1.

\begin{theorem}
The set of sections $\{f_T, T \in \mathcal{SS}(D,n) \}$
is a basis for $\Gamma(F_{\bk}(\C^n), \Lf)$.
\end{theorem}

\begin{example}\label{subdivide}

We consider the flag manifold  $F_\bk(\C^n) = F_{2,3}(\C^5)$ with the
very ample line bundle corresponding to the dominant weight
weight $2 \varpi_2 + \varpi_3$.
The associated Young diagram $D$ is

\begin{center}
D=
\begin{tabular}{| c | c | c |}
\hline
  &   &   \\ \hline
  &   &   \\ \hline
  \\ \cline{1-1}
\end{tabular}\ .
\end{center}

\vskip 12pt

Let $T$ be the filling of $D$ given by

\vskip 12pt

\begin{center}
T=
\begin{tabular}{| c | c | c |}
\hline
2  & 1  & 2  \\ \hline
3  & 4  & 5  \\ \hline
5  \\ \cline{1-1}
\end{tabular}.
\end{center}

Hence
\begin{center}  
$T_1 = $
\begin{tabular}{| c | c |}
\hline
1  & 2  \\ \hline
4  & 5  \\ \hline
\end{tabular}
\ and \ $T_2 = $
\begin{tabular}{| c |}
\hline
2 \\ \hline
3 \\ \hline
5 \\ \hline
\end{tabular}.  
\end{center}

\vskip 12pt

The  section $f_T$ is the pull-back to the flag manifold 
of the section over $Gr_2(\C^5) \times Gr_3(\C^5)$ given by 
$f_{T_1} \otimes f_{T_2}$, where
$f_{T_i}$ is the  section of $\mathcal{L}_{k_i}^{a_i}$ over
$Gr_{k_i}(\C^n)$.
\end{example}

\subsection{Duality of tableaux}

We define a map $\ast$ on tableaux as follows.  Let $T$ be a tableau.
Suppose the $i^{th}$ column of $T$ is $c_i = (p_1,...,p_\ell)$, with
distinct $p_j$'s.  Let $d_i = (q_1,...,q_{n -\ell})$ where 
$\{p_1,...,p_\ell,q_1,...,q_{n-\ell}\} = \{1,...,n\}$, and $q_t <
q_{t+1}$ for all $t$.  Let $*T$ be the tableau whose $i^{th}$ column
is $d_{n-i+1}$ for $1 \leq i \leq n$. 

\begin{example}

\begin{center}
$T=$
\begin{tabular}{| c | c | c |}
\hline
2  & 1  & 2  \\ \hline
3  & 4  & 5  \\ \hline
5  \\ \cline{1-1}
\end{tabular}
$\Longrightarrow$
\begin{tabular}{| c | c | c |}
\hline
    2  &     1  &     2  \\ \hline
    3  &     4  &     5  \\ \hline
    5  & \it{2} & \it{1} \\ \hline
\it{1} & \it{3} & \it{3} \\ \hline
\it{4} & \it{5} & \it{4} \\ \hline
\end{tabular}
$\Longrightarrow  *T = $
\begin{tabular}{| c | c | c |}
\hline
1  &  2  &  1   \\ \hline
3  &  3  &  4   \\ \hline
4  &  5  \\ \cline{1-2}
\end{tabular}

\end{center}
\end{example}

\begin{theorem}\label{semistandard_duality}
The map $\ast$ takes semistandard tableaux to semistandard tableaux.
\end{theorem}

\begin{proof}
Let $[n]$ denote the set $\{1,...,n\}$.  We define a 
partial order on the subsets of $[n]$ 
given by $I \leq J$ iff $|I| \geq |J|$ and $i_q \leq j_q$ for
all $q \leq t$, where $I = \{i_1,...,i_s\}$ and $J =
\{j_1,...,j_t\}$ have elements listed in strictly increasing order.
Define $*I$ as the complement of $I$ in $[n]$.  
We show that $I \leq J$ implies $*I \geq *J$ by induction on $n$.  If
$n=1$ this is trivial.  Now suppose that the statement is true for
$n-1$.  Let $I_{n-1} = I \cap [n-1]$ and let $J_{n-1} = J \cap
[n-1]$.  Define $*I_{n-1} = [n-1] \setminus I_{n-1}$ and $*J_{n-1} =
[n-1] \setminus J_{n-1}$.   Clearly $I_{n-1} \leq J_{n-1}$. 
By the induction hypothesis, $*I_{n-1} \geq
*J_{n-1}$.  Suppose there is some $*j_q > *i_q$, where $*j_q$ is the
$q^{th}$ element of $*J$ and $*i_q$ is the $q^{th}$ element of $*I$.  
Since $*J_{n-1} \leq *I_{n-1}$, we have that $|*J_{n-1}| \geq
|*I_{n-1}|$ and thus $q = 1+|*I_{n-1}|$ and
so $n = *i_q \in *I$, a contradiction with $*j_q > *i_q$.  

The columns of $*T$ are stictly increasing by definition.  Hence we
need only show that the rows are weakly increasing. 
Let $*t_{i,j}$ denote the $(i,j)$--entry of $*T$. We must show that
$*t_{i,j} \leq *t_{i,j+1}$ for all $(i,j)$ in the valid range.  
Let $d_j,d_{j+1}$ be adjacent columns in $\ast T$.  Then the
respective complementary columns $c_{n-j+1},c_{n-j}$ are 
adjacent columns of $T$.  Since $T$ is semistandard, the sets $I,J$ of
the entries of $c_{n-j},c_{n-j+1}$ respectively are such that 
$I \leq J$ for the
partial order on subsets of $[n]$ mentioned above.  Hence $*I \geq
*J$, and since $*I$ corresponds to $d_{j+1}$ and $*J$ corresponds to
$d_j$, we have that $*t_{i,j} \leq *t_{i,j+1}$ for all $i$.  
\end{proof}

Note that if the columns of $T$ are strictly increasing, then $\ast
\ast T = T$.

\vskip 12pt

We conclude our discussion of duality of tableaux
with a formula for how $\ast$ changes the weights. 
Let $w_i(T)$ be the number of times
the index $i$ appears in $T$, and let $wt(T) = (w_1(T),...,w_n(T))$.    
Let $D_{\ba}$ be the Young diagram with $m$ rows so that the $i$-th row
has length $a_1+a_2+\cdots + a_{m-i+1}$.

\begin{theorem}\label{tableauxweights}
For all tableaux $T$ with diagram $D_{\ba}$  
$$wt(T) + wt(\ast T) = \mathbf{\Lambda} = (|\ba|,|\ba|,...,|\ba|).$$ 
\end{theorem}

\begin{proof}
The $j^{th}$ column of $T$ and the $(n-j+1)^{th}$ column of $\ast T$
partition the set $\{1,...,n\}$.    The
total number of columns in either tableau is $|\ba| = \sum_i a_i$.
Fix any $i \in \{1,...,n\}$.
Let $c_j$ denote the $j^{th}$ column of $T$ and let $d_j$ denote the
$(n-j+1)^{th}$ column of $\ast T$.  The index $i$ is in exactly
one of $c_j, d_j$.  Hence $w_i(T) + w_i(\ast T) = |\ba|$, the total number
of columns.  
\end{proof}

\subsection{Duality of Mumford quotients for flag manifolds}

\subsubsection{The fundamental diagram}

The orthogonal complement map $\Psi$ maps the flag
$F_{k_1,...,k_m}(\C^n)$ to $F_{l_1,..,l_m}(\C^n$ where $l_i = n -
k_{p-i+1}$.  Denote $\bl = (l_1,...,l_m)$. 
Let $\mathbf{\Lambda} = (|\ba|,|\ba|,...,|\ba|)$ and
let $\bb = (a_m,a_{m-1},...,a_2,a_1)$, the tuple $\ba$ in reverse order, 

The extension of our duality theorem from Grassmannians to flag manifolds will 
result
from considering the diagram

\vskip 12pt

\begin{center}
\begin{math}
\begin{CD}
F_\bk(\C^n) @>{\Psi}>> F_\bl(\C^n) \\
@V{i}VV                  @VV{i}V \\
Gr_{k_1}(\C^n) \times \cdots \times Gr_{k_m}(\C^n)   
@>{F \circ \prod_i \Psi_i}>> 
Gr_{l_1}(\C^n) \times \cdots \times Gr_{l_m}(\C^n)
\end{CD}
\end{math}
\end{center} 

\vskip 12pt

We will refer to this diagram as the {\it fundamental diagram} in what follows.

\subsubsection{The bundle isomorphism $\hat{\Psi}$ and the induced 
isomorphism  of rings of sections}

\begin{lemma}
 There is a bundle isomorphism $\hat{\Psi}$  covering the isomorphism  
$\Psi:F_\bk(\C^n) \to F_\bl(\C^n)$. Moreover $\hat{\Psi}$ satisfies
$$\hat{\Psi} \circ h = h^{-1}\circ \hat{\Psi}.$$
\end{lemma}

\begin{proof}
We have seen in our analysis of duality for Grassmannians that the
isomorphism $\Psi_i$ can be covered by a bundle isomorphism 
$\hat{\Psi}_i$ satisfying
$$\hat{\Psi} \circ h = h^{-1}\circ \hat{\Psi}.$$ 

Hence the isomorphism $F \circ \Pi_i \Psi_i$ is
covered by the bundle isomorphism $F \circ \Pi_i \hat{\Psi}_i$
which satisfies the above equivariance condition with respect to
the product $H \times \cdots \times H$ and hence a fortiori with
respect to the diagonal.  But
by Lemma \ref{pullback} the bundles $\Lf$ and $\mathcal{L}_{\bk}^{\bb}$
are pull-backs by $i$ of the corresponding bundles on the products of
Grassmannians. Hence the pull-back of $F \circ \Pi_i \Psi_i$ by $i$ is a bundle 
isomorphism from $\Lf$ to $\mathcal{L}_{\bk}^{\ba}$.
\end{proof}

We obtain  induced isomorphisms $\widetilde{\Psi}:\Gamma(\Fk,{\Lf}^{\otimes N})
\to \Gamma(\Fl,{\Lg}^{\otimes N})$ by the formula
$$\widetilde{\Psi}(s)(x) = \hat{\Psi}(s(\Psi^{-1})).$$

Later we will need that $\widetilde{\Psi}(s)(x)$ intertwines
the representation $\rho_{\bk}$ with the action $\rho_{\bl}^{\theta}$
where $\rho_{\bl}^{\theta} = \rho_{\bl} \circ \theta.$ This follows
immediately from Lemma \ref{intertwining2}. We state this result
as a lemma.

\begin{lemma}\label{intertwining3}
Let $\rho_{N,a,\bk}$ be the representation of $GL_n(\C)$
on the vector space of sections 
$\Gamma(Gr_k(\C^n), \mathcal{L}_{\bk}^{\otimes aN})$.
Let $\rho_{N,a,\bk}^{\theta} = \rho_{N,a,\bk} \circ \theta$. 
Then we have
$$\widetilde{\Psi}_{N,a,\bk} \circ \rho_{N,a,\bk} \circ 
\widetilde{\Psi}_{N,a,\bl} = 
\rho_{N,a,\bl}^{\theta}.$$
\end{lemma}

We next compute the action of the  ring isomorphism $\widetilde{\Psi}$
on the elements $f_T$ and thereby determine how it changes weights of $H$.

\begin{theorem}\label{actionontableaux}

$$\widetilde{\Psi}(f_T) = \epsilon_T  f_{\ast T}.$$
\end{theorem}

\begin{proof} 
We have the following diagram of homogeneous coordinate rings corresponding to 
the 
fundamental diagram.

\begin{center}
\begin{math}
\begin{CD}
\Rk @>\widetilde{\Psi}>> \Rl \\
@A{i^*}AA                  @AA{i^*}A \\
\widetilde{\Rk} @>{ F \circ \prod_i \widetilde{\Psi}_i}>> \widetilde{\Rl}
\end{CD}
\end{math}
\end{center} 

Let $T$ be a tableau and $f_T \in \bar{R}_\bk$.  Then there are
tableaux $T_1,...,T_m$ such that 
$f_T = i^*(f_{T_1} \otimes \cdots \otimes f_{T_m})$, and
$\widetilde{\Psi}(f_T) = i^*(\widetilde{\Psi_1}(f_{T_1}) \otimes
\cdots \otimes \widetilde{\Psi_m}(f_{T_m})) = 
(\prod_i \epsilon_{T_i}) i^*(f_{(*T_1)} \otimes \cdots \otimes
f_{(*T_m)}) = \epsilon_T f_{*T}$.

\end{proof}
\begin{corollary}\label{weightchange}
$\widetilde{\Psi}$ maps the subspace of $\Rk^{(N)}$ of $H$--weight $N\br$ 
isomorphically to the subspace of $\Rl^{(N)}$ of 
$H$--weight $N(\mathbf{\Lambda} - \br)$.
\end{corollary}

We have now proved one of our main theorems.

\begin{theorem}
The isomorphism $\widetilde{\Psi}$ 
induces an isomorphism of graded rings

$$ \bigoplus_{N=0}^\infty 
\Gamma(F_\bk(\C^n),\mathcal{L}_\bk^{N \ba}(N \br))^H \cong 
\bigoplus_{N=0}^\infty 
\Gamma(F_\bl(\C^n),\mathcal{L}_\bl^{N \bb}(N \bs))^H$$

and consequently an isomorphism of Mumford quotients

$$F_\bk(\C^n)//_{\br} H \cong F_\bl(\C^n)//_{\bs}H.$$
\end{theorem}
\begin{proof}
The theorem follows immediately Corollary \ref{weightchange}.
\end{proof}

\subsection{Duality of K\"ahler structures}
We first observe that it follows from the fundamental diagram
on the corresponding result for Grassmannians that 
$\widetilde{\Psi}:\Fk \to \Fl$ is a holomorphic isometry. 
We now check that the induced map on quotient is also a holomorphic
isometry (hence symplectic). We already know it is holomorphic.
We have only to check that it is symplectic.

\begin{theorem}\label{flagsymplecticduality}
The map $\Psi$ induces a homeomorphism of the symplectic quotients:
$$\overline{\Psi }: F_\bk(\C^n) //_\br T \to F_\bl(\C^n)//_\bs T.$$
Furthermore, if $\br$ is a regular value of the momentum mapping, then
the symplectic quotients are symplectic manifolds and $\overline{\Psi}$ is a
symplectomorphism.  
\end{theorem}

\begin{proof}

We recall the fundamental diagram.
\begin{center}
\begin{math}
\begin{CD}
F_\bk(\C^n) @>{\Psi}>> F_\bl(\C^n) \\
@V{i}VV                  @VV{i}V \\
Gr_{k_1}(\C^n) \times \cdots \times Gr_{k_m}(\C^n)   
@>F \circ {\prod_i \Psi_i}>> 
Gr_{l_1}(\C^n) \times \cdots \times Gr_{l_m}(\C^n)
\end{CD}
\end{math}
\end{center} 

The map $\prod_i \Psi_i$ on the bottom is symplectic, and the inclusion
maps are symplectic, so the $\Psi$ map on the top is
symplectic.   The product $T^m$ acts in a Hamiltonian fashion on
each of the two products of Grassmannians.  Let $T^\Delta$ be the
diagonal.  Then the inclusion map $i$ is equivariant with respect to
$T^\Delta$.  Let $\mu^\Delta$ be the momentum mapping for $T^\Delta$.  
Hence it suffices to prove that $\prod_i \Psi_i$ carries 
$(\mu^\Delta)^{-1}(\br)$ to $(\mu^\Delta)^{-1}(|\ba|-\br)$.  
Choose $\mathbf{x} = (x_1,...,x_m) \in F_\bk$ such that
$\mu^\Delta(\mathbf{x}) = \br$.  Assume $\mu_i(x_i) = \br^{(i)} \in \R^n$.  
We
have $\mu^\Delta = \sum \mu_i$ where $\mu_i$ is the momentum mapping
for the action of $T_i$ on the $i^{th}$ factor of the product.  We
have seen that $\Psi_i$ takes $\mu_i^{-1}(\br^{(i)})$ to
$\mu_i^{-1}(\bs^{(i)})$ where $\br^{(i)} + \bs^{(i)} = (a_i,...,a_i)$.
Thus $\mu^\Delta((\prod_i \Psi_i)(\mathbf{x})) = (\sum_i a_i)
(1,1,...,1) - \br$. 
\end{proof}

\subsection{Duality of Gelfand-Tsetlin systems}\label{dualityofGTs}

\subsubsection{The duality map as a map of coadjoint orbits}\label{GTs}
In this subsection we will describe the duality map $\Psi$ as a map
of coadjoint orbits.  We identify the dual $\mathfrak{u}^{\ast}(n)$ with
the space $\mathcal{H}_n$ of $n$ by $n$ Hermitian matrices using the imaginary
part of the trace form. We note that $\varpi_n$ is identified to $I_n$.
We will in fact compute with the map $\Phi$
which operates of flags by taking 
the orthogonal complement relative to the positive definite 
Hermitian form $F$. We note that if we use $\sigma$ to denote
complex conjugation we have
$$\Phi = \Psi \circ \sigma.$$
The advantage in using $\Phi$ is that $\Phi$ is $U(n)$--equivariant.
 
 \smallskip
 
 Let $Flag$ denote the disjoint union of all the flag varieties
 of various lengths. We define a map $\mathcal{E}:\mathcal{H}_n \to Flag$
 as follows. Let $A \in \mathcal{H}_n$ be given. Let $\lambda_{i_1}, ..., 
 \lambda_{i_k}$ be the distinct eigenvalues of $A$ { \em arranged in decreasing
 order} and let $E_j$ be the eigenspace belonging to $\lambda_{i_j}$.
 Then we define $\mathcal{E}(A)$ to be the flag of partial sums
 of the $E_j$ whence
 $$\mathcal{E}(A) = (E_1, E_1 + E_2, ..., E_1+E_2+\cdots + E_{k-1}).$$
 
 It is important to note that $\mathcal{E}$ loses information. A  flag manifold equipped with an invariant 
symplectic form does not determine a unique orbit. If we change the orbit by adding
a multiple of $I_n$ we do not change the flag manifold as
a symplectic manifold. Indeed the map from the orbit to
the flag manifold assigns the flag attached to increasing partial sums of 
eigenspaces and the symplectic form depends only on the differences of
the eigenvalues. Note that $\mathbb{R}$ acts on $\mathcal{H}_n$ 
by translating by multiples of $I_n$. If we choose a cross-section to this action 
we can lift $\Phi$ and $\Psi$. We will henceforth choose the cross-section 
$\mathcal{H}_n^0$ of Hermitian matrices with smallest eigenvalue equal to zero. 

Let $\Xi: \mathcal{H}_n \to \mathcal{H}_n$ be the map given by
$$\Xi(A) = \lambda_{max}(A)I_n - A.$$
The reader will check that $\Xi$ is a Poisson map and carries
$\mathcal{H}_n^0$ into itself. Let $\Sigma:\mathcal{H}_n \to \mathcal{H}_n$
be complex conjugation so $\Sigma(A) = \overline{A}$.

  \smallskip
  
 \begin{lemma}
 The following diagram commutes
  
  \begin{center}
\begin{math}
\begin{CD}
\mathcal{H}_n^0 @>{\Xi}>> \mathcal{H}_n^0\\
@V{\mathcal{E}}VV                  @VV{\mathcal{E}}V \\
Flag   
@>{\Phi}>> 
Flag
\end{CD}
\end{math}
\end{center}

\end{lemma}  
  
\begin{proof}
Suppose that $\lambda_{i_1}, ..., 
 \lambda_{i_k}$ are the distinct eigenvalues of $A$ arranged in 
decreasing order. Then the eigenvalues of $-A$ arranged in
decreasing order are $- \lambda_{i_k}, ..., 
 -\lambda_{i_1}$ and consequently
 $\mathcal{E}(-A) = \Phi(\mathcal{E}(A)).$
 \end{proof}
 
\begin{corollary}
The following diagram commutes

  \begin{center}
\begin{math}
\begin{CD}
\mathcal{H}_n^0 @>{\Xi\circ \Sigma}>> \mathcal{H}_n^0\\
@V{\mathcal{E}}VV                  @VV{\mathcal{E}}V \\
Flag   
@>{\Psi}>> 
Flag
\end{CD}
\end{math}
\end{center}

\end{corollary}  

\begin{remark}
We note that (because $A^{\ast} =A$)
$$\Xi\circ \Sigma (A) = \lambda_{max}(A)I_n - A^t= \lambda_{max}(A) I_n +\theta(A).$$
Thus the duality map $\Psi$ lifted to the space of normalized coadjoint orbits
is once again given by the Chevalley involution (up to a translation).
\end{remark}

\subsubsection{Duality of Gelfand-Tsetlin systems}
We begin by recalling the definition of the Gelfand-Tsetlin 
Hamiltonians $\lambda_{i,j}: \mathcal{H}_n \to \mathbb{R}, 1 \leq i,j \leq n$. 
Let $B_j(A), 1 \leq j \leq n, $ be the upper principal $j$ by $j$ block of $A$. 
Then $\lambda_{i,j}(A)$ is the $i$--th eigenvalue of $B_j(A)$ (the eigenvalues
are arranged in (weakly) decreasing order). The $\lambda_{i,j}$'s
Poisson commute \cite{GuilleminSternberg} and the $\lambda_{i,n}$ are
Casimirs. The set of $\lambda_{i,j}$'s is called the Gelfand-Tsetlin system.
By restricting the Gelfand-Tsetlin system to any orbit we obtain an
integrable system on that orbit. Moreover the functions $\lambda_{i,j}$
descend to the torus symplectic quotients of orbits and hence define
integrable systems on the torus quotient of flag manifolds
$F_\bk(\C^n) //_\br T$.

\begin{theorem}
Assume that $F_\bk(\C^n)$ has the symplectic form corresponding to 
$a_1\varpi_{k_1} + \cdots + a_m \varpi_{k_m}$ and 
$F_\bl(\C^n)$ has the symplectic form corresponding to 
$b_1\varpi_{l_1} + \cdots + b_m \varpi_{l_m}$
Under the duality isomorphism
$\overline{\Psi} : F_\bk(\C^n) //_\br T \to F_\bl(\C^n)//_\bs T$
we have
$$\overline{\Psi}^{\ast}\lambda_{i,j}= |\bb| - \lambda_{j+1 - i,j} = 
|\ba| - \lambda_{i,j}.$$ 
\end{theorem}

\begin{proof}
It suffices to prove the above formula when $\lambda_{i,j}$ is pulled back to the 
orbit and so $\overline{\Psi}$ is replaced by
$\Xi$.  We have
$$\Xi^{\ast}\lambda_{i,j}(A) = \lambda_{i,j}(\Xi(A))= 
\lambda_i(B_j(\lambda_{max}(A)I_n) -A))= \lambda_i(|\ba|I_j - B_j(A)).$$
The theorem follows.
\end{proof}

\section{Self-duality}
Let $M$ be a flag manifold $G/P$ and $M^{opp} = G/P^{opp}$.
In this section we investigate the duality map 
$\overline{\Theta} :M//_{\br} H \to 
M^{opp}//_{\bs} H$ in case $P$ is conjugate to $P^{opp}$ and $\br = \bs$. 
Our main goal is to find the conditions when such a 
self-duality is trivial i.e. $\overline{\Theta} = Id$.
Roughly the following theorems say that self-duality is almost never trivial.
For our analysis of the case of $GL_n(\C)$ we will take advantage of
the solution of the quantum problem in \cite{MillsonToledano} 
although the analysis we give below would work for $GL_n(\C)$
as well. Our strategy below will be to first identify those cases with
integral $\ba$ and $\br$ for which duality is trivial and then show
that by scaling by real numbers we obtain all the real cases as well.

\subsection{The existence of good representations}

Recall that we have defined a dominant weight $\lambda$ (or representation
$V_{\lambda}$) to be good if $V_{\lambda}$ is self-dual and
if there exists $N$ such that the Chevalley involution does not
act as a scalar on $V_{N\lambda}[0]$. We will see shortly that this
condition on a representation exactly captures nontriviality
of the classical duality $\overline{\Theta}$ on the corresponding weight
variety. The point of this subsection is to prove that good representations
abound. We have

\begin{definition}\label{defCartanproduct}
Suppose $\lambda$ and $\mu$ are dominant weights and $V_{\lambda}$
and $V_{\mu}$ are the corresponding irreducible representations.
Then the Cartan product of $V_{\lambda}$ and $V_{\mu}$ is the
irreducible representation with highest weight $\lambda + \mu$.
There is a canonical surjection $\pi_{\lambda,\mu}:V_{\lambda} \otimes  V_{\mu} \to 
V_{\lambda + \mu}$. We will define the $N$--th Cartan power $C^N V_{\lambda}$
to be the irreducible representation $V_{N\lambda}$.
\end{definition}

We begin with a very useful lemma - the image of a nonzero 
{\em decomposable}
vector in the tensor product of two irreducibles in the Cartan product
is nonzero.

\begin{lemma} \label{decomposable}
Suppose that $\lambda$ and $\mu$ are dominant weights and
that $v_1 \in V_{\lambda}, v_2 \in V_{\mu}$ are nonzero vectors.
Then we have
$$\pi_{\lambda,\mu} (v_1 \otimes v_2) \neq 0$$
\end{lemma}
\begin{proof}
Use Borel-Weil to interpret $v_1$ and $v_2$ as sections $s_1$ and $s_2$ of line
bundles $\mathcal{L}_{\lambda}$ and $\mathcal{L}_{\mu}$ over
a flag manifold $M$.  The image 
$\pi_{\lambda,\mu} (v_1 \otimes v_2)$ then corresponds to the product 
$s_1 \cdot s_2$ of
the two sections under the multiplication map
$\Gamma(M, \mathcal{L}_{\lambda}) \otimes \Gamma(M, \mathcal{L}_{\mu})
\to \Gamma(M, \mathcal{L}_{\lambda} \otimes \mathcal{L}_{\mu})$.
But the product of two nonzero sections is never zero on an
irreducible variety.
\end{proof}

We first apply this to
\begin{lemma}\label{powers}
Suppose that $\theta$ does not act as a scalar on $V_{N_0\lambda}[0]$.
Then $\theta$ does not act as a scalar on $V_{k N_0\lambda}[0]$ for
any $k>0$.
\end{lemma}

\begin{proof}
Let $k>0$ be given. It will be convenient to argue in terms of sections.
By hypothesis there exists $s \in V_{N_0\lambda}[0]$ such that $s$
is not an eigenvector of $\theta$. We claim that $s^{\otimes k}$
is not an eigenvector of $\theta$ (note that since $M$ is
irreducible $s^{\otimes k}\neq 0$). Suppose to the contrary that
there exists $z$ with 
$$\theta(s^{\otimes k}) = z s^{\otimes k}.$$ 
Let $z_i, 1 \leq i \leq k$ be the $k$--th roots of $z$.
Then we have
$$\prod_{i=1}^k (\theta(s) -z_i s) = 0.$$
Again because $M$ is irreducible we must have
$\theta(s) -z_i s = 0$ for some $i$. This is a contradiction.
\end{proof}

We now prove that good representations are stable under Cartan product
and in fact much more is true.

\begin{theorem}\label{Cartanproduct}
Suppose that $V_{\lambda}$ and $V_{\mu}$ are self-dual representations
and $V_{\lambda}$ is good. Then the Cartan product $V_{\lambda+\mu}$
is good.
\end{theorem}

\begin{proof}
Since $V_{\lambda}$ is good there exists $N$ such that the Chevalley 
involution $\theta$ does not act on $V_{N\lambda}[0]$ as a scalar.
By the previous lemma $\theta$ does not act as a scalar on
$V_{kN\lambda}[0], k \geq 1$. Choose $k$ so that $V_{kN\mu}[0]\neq 0$.
Now choose a nonzero vector $v_+ \in V_{kN\lambda}[0]$ such that 
$\theta(v_+) = v_+$ and another nonzero vector $v_- \in V_{kN\lambda}[0]$
such that $\theta(v_-) = -v_-$. The Chevalley involution has either
$1$ or $-1$ as an eigenvalue on $V_{kN\mu}[0]$. For convenience assume
the former. Let $u$ be an eigenvector belonging to $1$. Then
the images of $v_+ \otimes u$ and $v_- \otimes u$ in the Cartan
product are nonzero by Lemma \ref{decomposable} and belong
to eigenvalues $1$ and $-1$ respectively.
\end{proof}

\subsection{The branching trick}

In this section we will give the main technique we will use below
to prove that certain fundamental representations are good. We will
refer to it as the ``branching trick''. It is (a refinement of) one of the
main techniques used in \cite{MillsonToledano}, see \S 4.3 and Proposition
4.6. 

We begin with the following lemma.

\begin{lemma}\label{branchingsemigroup}
Let $G_1$ and $G_2$ be simple complex Lie groups with $G_1 \subset G_2$. Let
$\lambda$ be a dominant weight for $G_2$ and $\mu$ be a dominant
weight for $G_1$.
If $V_{\mu}$ occurs in the restriction of the irreducible representation
$V_{\lambda}$ of $G_2$ to $G_1$ then $V_{N\mu}$ occurs in the restriction of 
the irreducible representation
$V_{N\lambda}$ of $G_2$ to $G_1$ for all $N \in \mathbb{N}$.
\end{lemma}

\begin{proof}
Suppose $v$ is a nonzero vector of weight $\mu$ in $V_{\lambda}$ which is
annihilated by the nilradical of $\mathfrak{g}_1$. Then the image
of the vector $v^{\otimes N}$ in $V_{N\lambda}$ is {\em nonzero}
(by Lemma \ref{decomposable}) has weight $N \mu$ and is
again annihilated by the nilradical of $\mathfrak{g}_1$.
\end{proof}

Recall that we have defined a  representation $V_{\lambda}$ to be good if 
$V_{\lambda}$ is 
self-dual and  for
some $N \in \mathbb{N}$ the Chevalley involution does not act
on $V_{N\lambda}[0]$ as a scalar. The ``branching trick'' is then
the following

\begin{proposition}\label{branchingtrick}
Suppose that $V_{\lambda}$ is an irreducible representation of simple complex Lie
group $G_2$ and that $G_1$ is a maximal rank subgroup so that the restriction
of $V_{\lambda}$ to $G_1$ contains an irreducible summand that is either
good or not self-dual. Then $V_{\lambda}$ is good.
\end{proposition}

\begin{proof}
Let $H$ be a Cartan subgroup of $G_2$ which is contained in $G_1$.
Let $\theta_1$ be a Chevalley involution of $G_1$ and $\theta_2$
be a Chevalley involution of $G_2$ such that both involutions
stabilize $H$ and consequently act on $H$ by inversion. Since $G_1$ is 
the centralizer of an element of $H$ it follows that $\theta_2$ carries $G_1$
into itself. Hence by Lemma \ref{characterizationofChevalley},
$\theta_1$ and $\theta_2$ are conjugate by an element $Ad h, h \in H$,
and by Lemma \ref{agreementofChevalley}, they coincide on the zero
weight space of any self-dual representation of $G_1$.

Suppose first that $V_{\mu}$ is good.
Since $V_{\mu}$ is good there exists $N$ so that $\theta_1$ does not act
as a scalar on $V_{N\mu}[0]$. Let $v \in V_{N\mu}[0]$ satisfy
$\theta_1(v) \neq \pm v$. Hence by the above paragraph
$\theta_2(v) \neq \pm v$. But by Lemma \ref{branchingsemigroup}, 
$V_{N\mu} \subset V_{\lambda}$ whence $V_{N\mu}[0] \subset V_{\lambda}[0]$
and $v \in V_{\lambda}[0]$.

Suppose now that $V_{\mu}$ is not self-dual. Choose $N$ such that
$N\mu$ is in the root lattice for $G_1$. Choose any nonzero vector $v \in V_{\mu}[0]$.
Then under the action of $\theta_1$ the vector $v$ goes to into a different
irreducible summand (corresponding to a copy of the dual of 
$V_{\mu}$) in the restriction of $V_{\lambda}$ to $G_1$.
This summand has intersection zero with $V_{\mu}$ by Schur's lemma.
Hence  once again we have $\theta_1(v) \neq \pm v$. The rest of
the argument is identical to that of the previous paragraph.
\end{proof}

\subsection{Quantum versus classical duality}
Our goal in this section is to compare the triviality of quantum
and classical self-dualities. 
Let $\lambda$ be a dominant weight
and $M$ be the corresponding flag manifold. We will identify the 
weight spaces
$V_{N\lambda}[N \br]$ with the $N$--th graded summand of the spaces of invariant
sections. In this section we will assume $\br=0$ and will 
identify the map  on sections $\widetilde{\Theta}$ with the action
of the Chevalley involution $\theta$ on the corresponding zero weight
space $V_{N\lambda}[0]$. We ask the reader to make the required modifications
in the proof to cover the case of $GL_n(\C)$ and the action of 
$\widetilde{\Psi}$
on the graded summands. Recall by Corollary \ref{weightchange} this
action corresponds with the affine involution 
$$\widetilde{\Psi}:V_{N\lambda}[\br] \to V_{N\lambda^{\vee}}[\mathbf{\Lambda} -\br].$$
We will see below that in this case the self-duality condition forces 
$\br = (|\ba|/2) \varpi_n = (1/2) \mathbf{\Lambda}$.

\begin{theorem}\label{quantumtoclassical}
Suppose that the symplectic manifold $M$ is self-dual and corresponds to an integral
orbit (the orbit of an element $\lambda = \lambda^{\vee}$ of the weight lattice).
Then the  classical duality map
$\overline{\Theta}$ is nontrivial (i.e.not equal to the identity) on 
$M//_0 H$ (resp. $M//_{\br}H$)
if and only if $V_{\lambda}$ is good.
\end{theorem}

The rest of this section will be devoted to proving the theorem.
One direction is easy. If $\theta$ acts as a scalar on every
summand then it is immediate that $\overline{\Theta}$ is equal
to the identity. Indeed 
we may choose $m$ such that the ring 
$S^{m} = \oplus_{k=0} ^{\infty} S^{(km)}$ is generated
by elements of degree one (i.e. $k=1$),see \cite{Bourbakicomalg},Chapter III,
\S 1.3, Proposition 3.  We rename
this ring $S$. By \cite{Dolgachev}, pg.39, we have an equality of
maximal projective spectra
$$Projm(S) = Projm(R).$$

Again $\theta$ acts by
a scalar on any graded summand of $S$ hence in
particular it acts as a scalar on the degree
one summand $V=S^{(1)}$.
Choose a basis of the degree one elements $V=S^{(1)}$ of $S$ to obtain a 
projective
embedding $F:M//H \to P(V^*)$.
We claim that $F$ satisfies
$$F(\overline{\Theta}(x)) = \theta(F(x)).$$
Indeed if we identify the dual of the projective space of the space of sections
with the hyperplanes in the space of sections then we have
$F(x) =H_x$ where $H_x$ is the hyperplane of sections 
that vanish at $x$,  see \cite{GriffithsHarris}, page 176.
The claim will follow if we show
$$H_{\overline{\Theta}(x)} = \widetilde{\Theta}^{-1}(H_x).$$
But $s(\overline{\Theta}(x)) = 0 
\Leftrightarrow \hat{\Theta}^{-1}(s(\overline{\Theta}(x)))
\Leftrightarrow \widetilde{\Theta}^{-1}(s)(x) =0$.
The claim follows.

By assumption $\theta$ acts as a scalar on $V$
and hence $\theta$ acts trivially on $P(V^*)$.
Since $F$ is injective we deduce from the equivariance formula
immediately above that $\overline{\Theta}$ is the identity.

\begin{remark} 
It is not enough to require that $\theta$ act as a scalar on the
degree one elements of the original graded ring $R$ because
$R$ might not be generated by elements of degree one.
\end{remark}

The rest of this section will be devoted to proving the converse i.e.
if there exists $N_0$ such that $\theta$ does not act as a scalar
on $V_{N_0\lambda}[0]$ then $\overline{\Theta}$ is 
not equal to the identity.

Accordingly we assume that  $\theta$ does not act as a scalar
on $V_{N_0\lambda}[0]$.
Replace the graded ring $R$ of $H$--invariant sections
by the subring $S$ given by $S = \oplus_{k=0}^{\infty} R^{kN_0}$.
By \cite{Dolgachev}, pg. 39, we have
an equality of maximal projective spectra
$$Projm(S) = Projm(R).$$
Also we note that by Lemma \ref{powers},  $\theta$ does not
act as a scalar on any graded summand of $S$. 

Finally, as before, we may choose $m$ such that the ring 
$S^{m} = \oplus_{k=0} ^{\infty} S^{(km)}$ is generated
by elements of degree one (i.e. $k=1$). We rename
this ring $S$. Again $\theta$ does not act as
a scalar on any graded summand of $S$ hence in
particular it does not act as a scalar on the degree
one summand and again by \cite{Dolgachev}, pg.39, we have 
$$Projm(S) = Projm(R).$$

Now we can complete the proof of the theorem.
Let $V$ be the space of degree one elements, $V=S^{(1)}$ of $S$. We obtain a 
projective,
embedding $F:M//_0 H \to P(V^{\ast})$ with $F(x) = H_x$ as above.
We have seen that  $F$ satisfies
$$F(\overline{\Theta}(x)) = \theta(F(x)).$$

Suppose now for the purpose of contradiction that $\overline{\Theta}$
is the identity. Then for all $x \in M//_0 H$ we have
$\theta(F(x)) = F(x)$. Let $<ImF>$ denote the smallest projective
subspace of $P(V^{\ast})$ containing the image of $F$.
We first check that $<ImF> = P(V^{\ast})$. Indeed, suppose that
$<ImF> \ \subsetneqq P(V^{\ast})$.  Then there exists
a nonzero element  $s \in V$ which pairs to zero with every element of
$<ImF>$. But this means that $\forall x \in M//_0 H, s(x) = 0$. This is a
contradiction.

Now we can prove that $\theta$ acts as a scalar on $V^{\ast}$ and hence on $V$.
Indeed we have the  eigenspace decomposition  
$V^{\ast} = (V^{\ast})^+ \oplus (V^{\ast})^-$. Hence
$P(V^{\ast})^+$ and $P(V^{\ast})^-$ are disjoint projective subspaces
of $P(V^{\ast})$ with union the fixed-point set of  $\theta$ on
$P(V^{\ast})$. Since $M$ is connected either
$ImF \subset P(V^{\ast})^+$ or $ImF \subset P(V^{\ast})^-$. 
Hence either $P(V^{\ast}) = <ImF> = P(V^{\ast})^+$ or 
$P(V^{\ast}) = <ImF> = P(V^{\ast})^-$. In either case
we find that $\theta$
acts as a scalar on $V$. 
This contradicts the assumption that $\theta$ does
not act as a scalar on any graded summand of $S$ and the theorem is proved.

\subsection{From general flag manifolds to Grassmannians}
In this section we prove the classical analogue of Theorem \ref{Cartanproduct}
that (in the case that
all representations of $G$ are self-dual)
will allow us to reduce to the study of $\overline{\Theta}$
from torus quotients of general flag manifolds  to ``Grassmannians'', that is
flag manifolds that are quotients $G/P$ where $P$ is {\em maximal}.
In case not all representations are self-dual the result will allow
us to reduce to the case of flag manifolds $G/P$ where $P$ is a
``next-to-maximal'' parabolic (see the treatment of $E_6$ below).

\smallskip

Let $P$ and $Q$ be parabolic subgroups of $G$ such that $P \subset Q$.
Then we have a quotient map $\pi:G/P \to G/Q$. Suppose that the flag
manifolds $G/P$ and $G/Q$ are carried into themselves by $\Theta$. 
We will abbreviate $G/P$ to $M$ and $G/Q$ to $N$.

\begin{lemma}
Suppose there exists $z \in N$ such that $\overline{H \cdot \Theta_N(z)} \cap 
\overline{H \cdot z} \neq \emptyset$. Then for any $w \in \pi^{-1}(z)$ we have
$$\overline{H\cdot\Theta_M(w)} \cap \overline{H\cdot w} \neq \emptyset.$$
\end{lemma}

\begin{proof} Suppose $\pi(w) =z$ and $x \in 
\overline{H\cdot\Theta_M(w)} \cap \overline{H\cdot w}$. Then because $\pi$
is a closed $H$--equivariant map we find $z \in \overline{H \cdot \Theta_N(z)}$,
a contradiction.
\end{proof}

Let $\overline{\Theta}_M$ resp. $\overline{\Theta}_N$ denote the induced
maps on the quotients by $H$ and $Fix(\overline{\Theta}_M)$ resp.
$Fix(\overline{\Theta}_N)$ denote their fixed-point sets. By the
previous lemma we have
$$ \pi^{-1}(Fix(\overline{\Theta}_N)) \subset Fix(\overline{\Theta}_M)).$$

We can now prove the reduction we need. 

\begin{proposition}\label{reductionproposition}
Suppose that $\Theta_N$ does not induce the identity on $N//_{\br}H$.
Then $\Theta_M$ does not induce the identity on $M//_{\br}H$.
\end{proposition}

\begin{proof}
Let $p:N^{ss} \to N//_{\br} H$ be the quotient map. Then, by assumption, 
$U = p^{-1}(Fix(\overline{\Theta_N})$ is a nonempty Zariski open subset
of $N$ whence $V =\pi^{-1}(U)$ is a nonempty Zariski open subset of $M$. 
Hence $V \cap M^{ss}$ is nonempty.  Let $w$ be a point in this intersection.
By the previous lemma we have 
$\overline{H\cdot\Theta_M(w)} \cap \overline{H\cdot w} \neq \emptyset$
and consequently the image of $w$ in $M//_{\br}$ is not fixed by 
$\overline{\Theta}_M$.
\end{proof}

\begin{remark}
The previous Proposition reduces the problem of showing that
$\overline{\Theta}$ is nontrivial on a torus quotient of a general flag 
variety $G/P$ to showing
that $\overline{\Theta}$ is nontrivial on the torus quotients of Grassmannians 
$G/Q$ where $Q$
is a maximal parabolic subgroup containing $P$.
\end{remark}

\subsection{Splitting the zero level - a nontriviality criterion}
In this subsection we give a useful condition we will use below,
see \S \ref{generalparameters}, to prove that 
$\overline{\Theta}$ is not equal to the identity for three special cases. 
We will assume that
$\theta$ is inner. Consequently $\Theta$ takes each Grassmannian
into itself.

\smallskip

Suppose that $f: M \to \prod_i ^n M_i$ is the inclusion from a flag
manifold into a product of Grassmannians and let $f_i:M \to M_i$ be the
surjection onto the $i$--th factor. Let $T \subset \prod_i^n T_i$
be the diagonal inclusion of the maximal compact torus into the product
of maximal compact tori and $\mu_i, 1\leq i \leq n$ the momentum map
for the action of $T_i$ on the $i$--th factor. Then $\mu = \sum_i ^n \mu_i$
is the momentum map for $T$ acting on the product. 

\begin{lemma}\label{splittingthezerolevel}
Suppose that $\theta$ is inner and there exists $x \in M$ such that 
\begin{enumerate}
\item $\mu(f(x)) = 0.$
\item For some $i, 1 \leq i \leq n$ we have $\mu_i(f_i(x)) = \br _i \neq 0.$
\end{enumerate}
Then the duality map $\overline{\Theta} : M//_0 H \to M//_0 H$  is not 
equal to the identity.
\end{lemma}

\begin{proof}
Observe that $\mu_i(f_i(\Theta(x))) = - \br_i$ and consequently
$f_i(\Theta(x))$ is not in the same $T_i$ orbit as $f_i(x)$. Hence
$x$ is not in the same $T$ orbit as $\Theta(x)$.
\end{proof}

\subsection{Self-duality for $SL_n(\C)$}
For those values of $\bk$ and $\br$
such that  $\bk = \bl$ and $\br = \bs$ the duality map $\Psi$ is a self-duality.
In this section we will prove that except for the case of $Gr_2(\mathbb{C}^4)$
with the symplectic form $2a \varpi_2$ and $\br = \varpi_4$ and one more
infinite family of examples (see below)the resulting self-duality maps are not
equal to the identity map. We first examine the condition 
$\br = \bs$. We assume that $F_\bk(\C^n)$ is equipped with the symplectic
form induced by embedding it as the orbit of $\sum_i a_i \varpi_i$.
Since $s_i = |\ba| -r_i$ the following formula is immediate.

\begin{lemma}
$\br = \bs \Rightarrow \br = (|\ba|/2) \varpi_n$
\end{lemma}

We need to know that $|\ba|/2$ is an integer. This is in fact
the case as will be seen in the following lemma.

Let $|\lambda|$ be the sum of the coefficients
of $\lambda$ when $\lambda$ is expressed in terms of the standard basis.
Recall that  $\lambda$ is in the root lattice if and only if $|\lambda|$ 
is divisible by $n$.
We now have

\begin{lemma}
Suppose $\lambda= \sum_{i=1}^{n-1} a_i \varpi_i$ is self-dual.
Then 
$$|\lambda|/n = |\ba|/2$$
and consequently if $\lambda$ is in the root lattice
then $|\ba| =\sum_{i=1}^{n-1} a_i$ is an even integer. 
\end{lemma}

\begin{proof}
Suppose first that  $n$ is odd, $n = 2m+1$. Since $a_i = a_{n-i}$ we
have 
$$|\lambda|= \sum_{i=1}^m a_i i + a_{n-i} (n-i) = (\sum_{i=1}^m a_i i) n
=(|\ba|/2) n.$$

Assume now  that $n=2m$. Then as in the odd case we have
$$ |\lambda| = 2m(\sum_{i=1}^{m-1}a_i) + ma_m = 2m(|\ba|/2)$$
\end{proof}

We now recall Theorem 7.2 of \cite{MillsonToledano}.

\begin{theorem}
Suppose $\lambda$ is a dominant weight for an irreducible representation
of $SL_n(\mathbb{C})$ which is in the root lattice. Then the Chevalley 
involution $\theta$ of $SL_n(\mathbb{C})$ acts as a scalar on
the zero weight space $V_{\lambda}[0]$ if and only if $\lambda$
is one of the following
\begin{enumerate}
\item $\lambda=(a,0,\ldots,0,-a)$, $a\in\mathbb{N}$.
\item $\lambda=(\underbrace{1,\ldots,1}_{k},0,\ldots,0,
\underbrace{-1,\ldots,-1}_{k})$, $0\leq k\leq n/2$.
\item $\lambda=(a,a,-a,-a)$, $a\in\mathbb{N}$.
\end{enumerate}
\end{theorem}

We want to deduce from this theorem the analogous result for the
action of $\widetilde{\Psi}$ on the graded summands of the
space of $H$--invariant sections which we know corresponds to the weight
space $V_{N\lambda}[(N/2) |\ba| \varpi_n]$. We refer the reader
to \S \ref{actionofChevalley} for the definition of the
action of the Chevalley involution on $V_{\lambda}$ and $V_{\lambda}[0]$
and the definition of the operator $\Theta_{V_{\lambda}}$.
\smallskip

We now show we may choose $\widetilde{\Psi}$ for $\Theta_{V_{\lambda}}$.
This is a consequence of the next lemma which in turn is an immediate
consequence of Lemma \ref{intertwining3}.

\begin{lemma}\label{intertwining4}
Let $\rho_{\lambda}$ be the representation of $GL_n(\C)$ on
$V_{\lambda}$. Assume that $V_{\lambda}$ is self-dual.
Then we have as representations of $SL_n(\C)$
$$\widetilde{\Psi} \circ \rho_{\lambda} \circ \widetilde{\Psi}
= \rho_{\lambda}\circ \theta$$
or 
$$\widetilde{\Psi} = \Theta_{V_{\lambda}}.$$
\end{lemma}

We next need to further modify the above theorem because.
we are normalizing the highest weight $\lambda$ of an irreducible
representation of $SL_n(\C)$ to have last component zero 
rather than to have
the sum of its components $|\lambda|$ equal to zero. Let $H_1$
denote the subgroup of $H$ of elements of determinant $1$. We note that $\lambda$ is 
in the root lattice if and only if $|\lambda|$ is
divisible by $n$ and then the  zero weight space for $H_1$ in $V_{\lambda}$ coincides
with the $H$ weight space $V_{\lambda}[(|\lambda|/n)\varpi_n]$
=$V_{\lambda}[(|\ba|/2)\varpi_n]$. We now obtain the version of the previous 
theorem that we need

\begin{corollary}
Suppose $\lambda$ is a dominant weight for an irreducible representation
of $GL_n(\C)$ which is self-dual and in the root lattice. Let $M$ be
the flag manifold corresponding to $\lambda$ and $\mathcal{L}$
be the line bundle over $M$ corresponding to $\lambda$. Then the action
of $\widetilde{\Psi}$ on the graded component of the graded ring of
$H$--invariant sections of the line bundle $\mathcal{L}$ corresponding to the vector space
$V_{\lambda}[\br]$ is a scalar if and only if either $n=2$ or 
$\lambda$ and $\br$ are one of
the following
\begin{enumerate}
\item $\lambda=a\varpi_1 + a\varpi_{n-1}$, $\br = a \varpi_n$,$a\in\mathbb{N}$.
\item $\lambda=\varpi_k + \varpi_{n-k}, 2 \leq k \leq n-2$, $\br =\varpi_n$
\item $n=4$ and $\lambda=2a\varpi_2$, $\br= a\varpi_4$, $a\in\mathbb{N}$.
\end{enumerate}
\end{corollary}

We now prove
\begin{theorem}
Assume that $\bk$ and $\br$ satisfy the self-duality conditions $\bk=\bl$
and $\br = \bs$.
The  self-duality $\overline{\Psi} : F_\bk(\C^n) //_\br H \to F_\bk(\C^n)//_\br 
H$ is 
equal to the identity
if and only if either $n=2$ or the flag manifold is 
\begin{enumerate}
\item  $F_\bk(\C^n) = F_{1,n-1}(\mathbb{C}^n)$ with the
symplectic form $a \varpi_1 + a \varpi_{n-1}$ and $\br = a \varpi_n$.
\item $F_\bk(\C^n) = Gr_2(\mathbb{C}^4)$ with the
symplectic form $2a\varpi_2 $ and $\br = a \varpi_4$.
 
\end{enumerate}
\end{theorem}
\begin{proof}
The theorem follows from Theorem \ref{quantumtoclassical}.
We note that since, by the above theorem, $\widetilde{\Psi}$ acts as a scalar on 
all graded summands
for the two exceptional cases we may apply (the easy direction of) Theorem 
\ref{quantumtoclassical}
to deduce that $\overline{\Psi}$ is in fact equal to the identity  in 
these two cases.
\end{proof}

\subsection{Self-duality for the isometry groups of bilinear forms.}
In this section we assume that $M$ is a flag manifold of a
classical group $G$ where $G$ is either a symplectic group or a
special orthogonal group. In the case of symplectic groups and odd orthogonal 
groups
$-1$ is an element of the Weyl group and all
dualities are self-dualities.  We will
use the theory of admissible pairs \cite{Lakshmibai} or \cite{Musili} to
construct certain basis elements in the graded summands. However our
goal is to find weight zero monomials that are not eigenvectors of $\theta$
and as we will see below in order to do this 
we need only those
standard monomials which are products of extremal weight vectors, i.e.
correspond to a {\em trivial} admissible pairs. Hence,we do 
not need the difficult part of
the theory which constructs elements $p(\lambda,\phi)$
for a nontrivial admissible pair $\lambda,\phi$. All we need
from the general theory is that the standard monomials formed
from Bruhat chains of extremal weight vectors are linearly independent.
We will also use that for the Grassmannians associated to $Sp_{2n}(\C)$
and $SO_{2n+1}(\C)$ the Bruhat order on the relative Weyl poset
$W^P$ coincides with the restriction of the Bruhat order from
the ambient linear group \cite{Lakshmibai}, page 363, IX and page 365, IX.
We will also need the corresponding fact for the Bruhat order on $W^P$
for $P$ the subgroup of $SO_{2n}(\C)$ which stabilizes one of the two
types of Lagrangian subspaces (so $P$ is miniscule), see 
\cite{GonciuleaLakshmibai}, pg. 158.

Our goal is to prove the following theorem.

\begin{theorem}\label{nontrivialitytheorem}

  \hfill

\begin{enumerate} 
\item Suppose $G = Sp_{2n}(\C)$.Then
the duality map $\overline{\Theta}$ is not equal to the identity
with the exception of the torus quotients of 
\begin{enumerate}
\item The projective space $\mathbb{CP}^{2n-1}$.
\item The Lagrangian Grassmannian $Gr_{2}^{0}(\C^4)$.
\end{enumerate}
In both (a) and (b) the map $\overline{\Theta}$ is the identity.
\item Suppose now that $G= SO_{2n+1}(\C)$. Then the duality map 
$\overline{\Theta}$  is not equal to the identity
with the exception of the torus quotients of 
\begin{enumerate}
\item The quadric hypersurface $\mathcal{Q}\subset \mathbb{CP}^{2n}$.
\item The Lagrangian Grassmannian $Gr_{2}^{0}(\C^5)$.
\end{enumerate}
In both (a) and (b) the map $\overline{\Theta}$ is the identity.
\item Suppose now that $G= SO_{2n}(\C)$. Then the duality map 
$\overline{\Theta}$ 
 is not equal to the identity
with the exception of the torus quotients of 
\begin{enumerate}
\item The quadric hypersurface $\mathcal{Q}\subset \mathbb{CP}^{2n-1}$.
\item The isotropic Grassmannian $Gr_2^0(\C^6)$.
\item The Lagrangian Grassmannians $Gr_{2}^{0}(\C^4)^{+}, 
Gr_{2}^{0}(\C^4)^{-},Gr_{4}^{0}(\C^8)^{+}, Gr_{4}^{0}(\C^8)^{-}$
\item The isotropic flag manifold $F_{1,2}^0(\C^{4})$. 
\end{enumerate}
In (a), (b) (c) and (d) the map $\overline{\Theta}$ is the identity.
\end{enumerate}

\end{theorem}

\section{Proof of theorem \ref{nontrivialitytheorem}}

This section will be devoted to proving the theorem. As explained about we will use the
theory of standard monomials (in tableaux form) due to Seshadri and Lakshmibai.

\smallskip

We will begin with a lemma that reduces the proof of the theorem for the
classical groups to the case of 
the Lagrangian Grassmannians or equivalently to the problem of 
when the last fundamental
representation is good. 

\begin{lemma}\label{squarecase}
If the $k$--th fundamental representation is good or not self-dual for $G = Sp_{2k}(\C)$, 
resp. $SO_{2k+1}(\C)$,
resp. $SO_{2k}(\C)$ then the $k$--th fundamental representation is good for $G = Sp_{2n}(\C)$, 
resp. $SO_{2n+1}(\C)$, resp. $SO_{2n}(\C)$ with $n >k$.
\end{lemma}

\begin{proof}
First observe that for $Sp_{2n}(\C)$ the fundamental representations
are the primitive exterior representations $\bigwedge_0^k(\C^{2n})$
 of the standard representation and for the orthogonal groups they
 are either the exterior powers of the standard representation or
 a $Spin$ representation. In this later case the Cartan square is
 an exterior power of the standard representation (or else contained
 in it as the subspace fixed by the complex Hodge star). 
 
 We claim we  have following formulas for branching to
 the maximal rank subgroups $Sp_{2k}(\C) \times Sp_{2(n-k)}(\C)$,
 resp. $SO_{2k+1}(\C) \times SO_{2(n-k)}(\C)$, resp. 
 $SO_{2k}(\C) \times SO_{2(n-k)}(\C)$.
  
\smallskip

\begin{enumerate}
\item $\bigwedge^k_0(\C^{2n})|Sp_{2k}(\C) \times Sp_{2(n-k)}(\C)$ \ contains 
\ $\bigwedge_0^k(\C^{2k}) \boxtimes \bigwedge^0(\C^{2n-2k})$
as representations of $Sp_{2k}(\C) \times Sp_{n-k}(\C)$.
\item $\bigwedge^k(\C^{2n+1})|SO_{2k+1}(\C) \times SO_{2(n-k)}(\C)$ 
\ contains \ 
$\bigwedge^k(\C^{2k+1}) \boxtimes \bigwedge^0(\C^{2n-2k})$
as representations of $SO_{2k+1}(\C) \times SO_{2(n-k)}(\C)$.
\item $\bigwedge^k(\C^{2n})|SO_{2k}(\C) \times
SO_{2(n-k)}(\C) $ \ contains \ both \  
$\bigwedge^k(\C^{2k})_{+} \boxtimes \bigwedge^0(\C^{2n-2k})$
 and
$\bigwedge^k(\C^{2k})_{-} \boxtimes \bigwedge^0(\C^{2n-2k})$
as representations of $SO_{2k}(\C) \times SO_{2(n-k)}(\C)$. 
\end{enumerate}

\smallskip

These formulas in turn are an immediate consequence of the
formula that the standard representation restricts to the direct
sum of the two representations obtained from the standard representation from 
one factor tensored with the
trivial representation from the other factor together
with the usual formula for the $k$--th exterior power of a direct sum.
Indeed to prove (1) we first note that the analogue holds for the full
exterior power $\bigwedge^k(\C^{2n})$. But in general we have
$\bigwedge^k_0(\C^{2n}) = 
\oplus_{j=0} ^k \bigwedge^j_0(\C^{2i}) \otimes \bigwedge^{k-j}_0(\C^{(2k-2i)})$
and the first formula follows.

The second formula is follows immediately from the usual formula for the exterior power of a 
direct sum. The third formula follows from the formula for restricting
an exterior power together with the fact that $\bigwedge^k(\C^{2k})$
is the direct sum of $\bigwedge^k(\C^{2k})_{+}$ and $\bigwedge^k(\C^{2k})_{-}$.

The lemma now follows from Proposition \ref{branchingtrick} and the 
observation  $C^N(U \boxtimes V) = C^N(U) \boxtimes C^N(V)$ where $U\boxtimes V$
is the outer tensor product of the irreducible representations $U$
of $G_1$ and $V$ of $G_2$ (here $C^N(W)$ denotes the $N$--th Cartan power
of $W$). 
\end{proof}

\subsection{The second fundamental representation for the symplectic
and  orthogonal groups}

Before beginning our study of the last fundamental representation(s)
for the classical groups we deal
with the  case of the Grassmannians $Gr_2^0(\C^n)$.

\begin{lemma} Let $V_{\varpi_2}$ denote the irreducible representation
of either $Sp_{2n}(\C)$ or $SO_{2n+1}(\C)$ 
with highest weight $\varpi_2$. 
Then 
\begin{enumerate}
\item $V_{\varpi_2}$ is good
for $Sp_{2n}(\C)$ provided $n \geq 3$.
\item $V_{\varpi_2}$ is good
for $SO_{2n+1}(\C)$ provided $n \geq 3$.
\item $V_{\varpi_2}$ is good
for $SO_{2n}(\C)$ provided $n \geq 4$. 
\end{enumerate}
Moreover in all three cases the Chevalley involution does not
act as a scalar on $V_{3\varpi_2}[0]$.
\end{lemma}

\begin{proof}
We first give a proof using standard monomials valid for the symplectic
and odd orthogonal cases.
The weight zero cubic  monomial $m_I$ in the Pl\"ucker coordinates given by
$m_I = X_{12} X_{3\overline{2}} X_{\overline{3}\overline{1}}$ is not an 
eigenvector of $J$. Hence
$J$ has eigenvalues of both signs on the third graded summand of
the homogeneous coordinate ring of $Gr_{2}^0(\C^{n})//_0 H$.

We now deal with the case of $SO_{2n}(\C)$. First we treat the case
of $SO_8(\C)$ by branching the third Cartan power to the maximal subgroup of 
maximal rank $(SL_2(\C))^4$. Each zero weight summand in the restriction
is the Cartan product $S^{k_1}(\C^2) \boxtimes S^{k_2}(\C^2) \boxtimes
S^{k_3}(\C^2) \boxtimes S^{k_4}(\C^2)$. Here $S^k(\C^2)$ denotes the
$k$--th symmetric power of the standard representation. The zero
weight space of this summand is nonzero if and only if each $k_i$
is even. The Chevalley involution acts on this summand as the outer tensor
product of the Chevalley involution in each factor and hence each
summand is an eigenvector of Chevalley (because the zero weight
space of each summand is one dimensional). The sign on the $i$--th
factor depends whether $k_i$ is congruent to $0$ or $2$ modulo $4$ ( it is $+1$
in the first case and $-1$ in the second). Using LiE we find
the $-1$ eigenspace is large and the $+1$ eigenspace is one-dimensional.

To do the case of general $n$ we branch to the maximal subgroup of maximal
rank $SO_{2n-2}(\C) \times SO_2(\C)$. By the exterior product of
a direct sum formula we have  
$\bigwedge^2(\C^{2n})|SO_{2n-2}(\C) \times SO_2(\C)$ contains
$\bigwedge^2(\C^{2n-2}) \boxtimes \bigwedge^0(\C^2)$ as a summand.
But this summand is good by induction.
\end{proof}

\begin{remark}
The representation $V_{\varpi_2}$ is not good for $SO_6(\C)$. Indeed the second
exterior power $\bigwedge^2(\C^6)$ is not a fundamental representation.
It is the Cartan product of the two spin representations,
it has highest weight $(1,1,0) = \varpi_2 + \varpi_3$. Under the
isomorphism between $Spin_6(\C)$ and $SL_4(\C)$ it corresponds to
the sum of the first and third fundamental representations which
in turn corresponds to the exceptional case of the flag manifold 
$F_{1,3}(\C^4)$. Thus the Chevalley involution acts as a scalar
on the zero weight spaces of all Cartan powers of $\bigwedge^2(\C^6)$.
\end{remark}

\subsection{The symplectic group $Sp_{2n}(\C)$ }
\hfill

In this subsection we will prove
\begin{lemma}
  The representation of $Sp_{2n}(\C)$ with highest weight $\varpi_n$
 is good if $n \geq 3$.
 
 \end{lemma}
 
 Note  that $J(e_i) = \epsilon e_{\overline{i}}$ where $\overline{i}= 2n+1 -i$
 and $\epsilon = -1$ if $1 \leq i \leq n$ and $+1$ otherwise. 
 Each column of the tableaux below represents the wedge of the coordinate
 vectors corresponding to the entries in that column. According $J$
 acts on the basis vector represented by the tableau by changing
 each entry to its bar and multiplying by a sign which will not be
 important to us here.

 \medskip
 The rest of this subsection will be devoted to proving the lemma. 
 
Consider the case of $Gr_3^0(\C^6)$, the space of isotropic three
dimensional subspaces of $\C^6$.  Let

\vskip 12pt
\begin{center}
$\alpha = $
\begin{tabular}{| c | c | c | c |}
\hline 1  & 1  & 2 & 3  \\ \hline 2  & 4  & 4 & 5  \\ \hline 3  &
5  & 6 & 6  \\ \hline
\end{tabular} $\quad$
$\beta = $
\begin{tabular}{| c | c | c | c |}
\hline 1  & 1  & 2 & 4  \\ \hline 2  & 3  & 3 & 5  \\ \hline 4  &
5  & 6 & 6  \\ \hline
\end{tabular} 
\end{center}

The sections $\alpha$ and $\beta$ are standard basis vectors (i.e.
correspond to standard monomials) of the irreducible representation
$V_{4 \varpi_3}$  since the 
columns correspond to {\em extremal} weights of the third exterior power of the
standard representation ( i.e. they index
isotropic coordinate planes) and they are increasing in the Bruhat
order. They
are in the $0$--weight space $V_{4 \varpi_3}[0]$, since the
indices 1 through 6 appear exactly twice each (more generally, an
index $i$ needs to appear with the same frequency as its
complement $\overline{i} = 2n + 1 - i$.) The Chevalley involution
$\theta$ maps $\alpha$ to $\beta$. Since $\alpha$ and $\beta$ represent
standard monomials, the sections corresponding to
$\alpha$ and $\beta$ are linearly independent.
Hence $\theta$ does not act as a scalar on $V_{4\varpi_3}[0]$.

We now construct analogous sections $\alpha_n$ and $\beta_n$
of the irreducible representation $V_{4 \varpi_n}[0]$ of
$Sp_{2n}(\C)$ for all $n \geq 3$ by
induction. Suppose that $n \geq 4$ and $\alpha_{n-1}, \beta_{n-1},
\gamma_{n-1}$ have already been constructed. Let $\alpha_n$ have
an $n$ by $4$ diagram, and let $\alpha_n(1,1) = \alpha_n(1,2) =
1$, $\alpha_n(n,3)=\alpha_n(n,4)=2n$.  For $j=1,2$, and $i \geq 2$
let $\alpha_n(i,j) = \alpha_{n-1}(i-1,j)+1$, and for $j=3,4$ and
$i \leq n-1$ let $\alpha_n(i,j) = \alpha_{n-1}(i,j)+1$.  In other
words, to get $\alpha_n$, first add $1$ to all the entries of
$\alpha_{n-1}$. Then slide the first two columns down one level,
and put in two $1$'s and two $2n$'s in the remaining empty slots.
Here is an example to get $\alpha_4$ from $\alpha_3$.

\begin{center}
$\alpha_3 =$
\begin{tabular}{| c | c | c | c |}
\hline 1  & 1  & 2 & 3  \\ \hline 2  & 4  & 4 & 5  \\ \hline 3  &
5  & 6 & 6  \\ \hline
\end{tabular}
$\longrightarrow$
\begin{tabular}{| c | c | c | c |}
\hline 2  & 2  & 3 & 4  \\ \hline 3  & 5  & 5 & 6  \\ \hline 4  &
6  & 7 & 7  \\ \hline
\end{tabular}
$\longrightarrow$
\begin{tabular}{| c | c | c | c |}
\hline
   &    & 3 & 4  \\ \hline
2  & 2  & 5 & 6  \\ \hline 3  & 5  & 7 & 7  \\ \hline 4  & 6  &
&    \\ \hline
\end{tabular}
$\longrightarrow$
\begin{tabular}{| c | c | c | c |}
\hline 1  & 1  & 3 & 4  \\ \hline 2  & 2  & 5 & 6  \\ \hline 3  &
5  & 7 & 7  \\ \hline 4  & 6  & 8 & 8  \\ \hline
\end{tabular}
$ = \alpha_4$
\end{center}
Since the columns of $\alpha_{n-1}$ are isotropic,  no
column of $\alpha_{n-1}$ contains both an index $i$ and its
complement $\overline{i} = (2n-1)-i$. Note that $\overline{i+1} = 2n+1 - (i+1) = ((2n-1) - i)
+ 1 = \overline{i} + 1$.  Therefore the columns of $\alpha_n$ are isotropic
as well. Sliding the first two columns down preserves the property
that rows are weakly increasing, and the weight of $\alpha_n$ is
$0$ since all indices occur twice.

The section $\beta_n$ is constructed from
$\beta_{n-1}$ in the same way as $\alpha_n$ is
constructed from $\alpha_{n-1}$. The Chevalley involution takes
$\alpha_n$ to $\beta_n$. Since $\beta_n$ is standard it is independent
of $\alpha_n$ and consequently the Chevalley involution has 
eigenvalues $+1$ and $-1$ on degree $4$ weight $0$ sections.
The lemma is now proved.

We have now shown that $\overline{\Theta}$ is nontrivial on the torus
quotients of the symplectic
Grassmannians $Gr_n^0(C^{2m})//_0 H$ for all $n \geq 3$.

\subsection{The orthogonal group $SO_{2n+1}(\C)$}

\hfill

In this subsection we will prove
\begin{lemma}
  The representation of $SO_{2n+1}(\C)$ with highest weight $\varpi_n$
 is good if $n \geq 3$.

 \end{lemma}

The construction is  similar to that of the symplectic case,
as one might expect since the Weyl groups are the same. The
difference is that the middle index $n+1$ may not appear in the
tableau.  To get $\alpha_n$ for the orthogonal group from the
$\alpha_n$ for the symplectic group, simply add one to each index
which is greater than or equal to $n+1$.  For example, when $n=3$,

\begin{center}
$\alpha_3 = $
\begin{tabular}{| c | c | c | c |} \hline
1 & 1  & 2 & 3  \\ \hline 2 & 5  & 5 & 6  \\ \hline 3 & 6 & 7 & 7
\\ \hline
\end{tabular} $\quad$
$\beta_3 = $
\begin{tabular}{| c | c | c | c |}
\hline 1  & 1  & 2 & 5  \\ \hline 2  & 3  & 3 & 6  \\ \hline 5  &
6 & 7 & 7  \\ \hline
\end{tabular}
\end{center}

\subsection{The orthogonal groups $SO_{4n}(\C)$}

\hfill

In this subsection we will prove
\begin{lemma}\label{orthogonalsquarecase}
 The irreducible representations of $SO_{4n}(\C)$ with highest weights
 $\varpi_{2n-1}$ or  $\varpi_{2n}$
 are good if  $n \geq 3$.
 
 \end{lemma}
  
  \medskip
  
  The Grassmannian of isotropic $n$-dimensional spaces in $\C^{2n}$
has two components, $Gr_n^0(\C^{2n})^+$ and $Gr_n^0(\C^{2n})^-$.
The corresponding representations are the two spin representations 
$\Delta_{2n}^+$ and $\Delta_{2n}^-$ and are miniscule.
The weights of each representation lie in a single Weyl group orbit
and consequently must have the same
parity in the number of negative signs, since any Weyl group
element must negate an even number of components. Without
loss of generality we will treat the case of $Gr_n^0(\C^{2n})^+$.
The (extremal) standard monomials correspond to
tableaux that have columns representing weights which are all in the
same Weyl orbit, and increasing in the Bruhat order induced from 
$SL_{2n}(\C)$, see \cite{GonciuleaLakshmibai}, page 158, for a
description of the Bruhat poset.

We now show that  the Chevalley involution does not act as a scalar for 
$n$ even, $n \geq 6$. To prove this we begin with $Gr_6^0(\C^{12})$.  
Let 

\vskip 12pt

\begin{center}
$\alpha_6 = $
\begin{tabular}{| c | c | c | c |}
\hline 
1  &  1 &  2 &  5  \\ \hline 
2  &  3 &  4 &  6  \\ \hline 
3  &  7 &  7 &  9  \\ \hline 
4  &  8 &  8 & 10  \\ \hline 
5  &  9 & 10 & 11  \\ \hline 
6  & 11 & 12 & 12  \\ \hline 
\end{tabular} $\quad$
$\beta_6 = $
\begin{tabular}{| c | c | c | c |}
\hline 
1  &  1 &  2 &  7  \\ \hline 
2  &  3 &  4 &  8  \\ \hline 
3  &  5 &  5 &  9  \\ \hline 
4  &  6 &  6 & 10  \\ \hline 
7  &  9 & 10 & 11  \\ \hline 
8  & 11 & 12 & 12  \\ \hline 
\end{tabular}  $\quad$
\end{center}

The reader will observe that $\alpha_6$ and $\beta_6$
satisfy that each column has even parity, and so they define 
standard basis elements of the fourth Cartan power of the even 
spin representation $\Delta^+_{12}$ of $Spin(12)$. Furthermore, 
$\alpha_6$ is mapped by $\theta$ to $\beta_6 \neq \alpha_6$.

We construct $\alpha_{2k}$ and $\beta_{2k}$ for $k \geq 3$.
To get $\alpha_{2k}$ from $\alpha_{2k-2}$, first add $2$ to each 
entry, then slide the first two columns down two levels, and put 
in two each of $1,2,4k-1,4k$ in the appropriate positions.

For example,
\begin{center}
$\alpha_6 = $
\begin{tabular}{| c | c | c | c |}
\hline 
1  &  1 &  2 &  5  \\ \hline 
2  &  3 &  4 &  6  \\ \hline 
3  &  7 &  7 &  9  \\ \hline 
4  &  8 &  8 & 10  \\ \hline 
5  &  9 & 10 & 11  \\ \hline 
6  & 11 & 12 & 12  \\ \hline 
\end{tabular}
$\longrightarrow$
\begin{tabular}{| c | c | c | c |}
\hline 
   &    &  4 &  7  \\ \hline 
   &    &  6 &  8  \\ \hline 
3  &  3 &  9 & 11  \\ \hline 
4  &  5 & 10 & 12  \\ \hline 
5  &  9 & 12 & 13  \\ \hline 
6  & 10 & 14 & 14  \\ \hline 
7  & 11 &    &     \\ \hline
8  & 13 &    &     \\ \hline
\end{tabular} 
$\longrightarrow$
\begin{tabular}{| c | c | c | c |}
\hline 
1  &  1 &  4 &  7  \\ \hline 
2  &  2 &  6 &  8  \\ \hline 
3  &  3 &  9 & 11  \\ \hline 
4  &  5 & 10 & 12  \\ \hline 
5  &  9 & 12 & 13  \\ \hline 
6  & 10 & 14 & 14  \\ \hline 
7  & 11 & 15 & 15  \\ \hline
8  & 13 & 16 & 16  \\ \hline
\end{tabular} 
$ = \alpha_8$  
\end{center}

The parity is still even since adding $2$ does not change the parity 
as you go from $2k-2$ to $2k$, and adding the indices $1,2$ does not 
affect the number of negative signs of the weights associated to
the first two columns.  Adding $4k-1,4k$ to the last two columns adds
two negative signs, and thus parity remains even.

The section $\beta_{2k}$  is formed in the same manner from the
 $\beta_{2k-2}$.  Hence the Chevalley involution is
non-trivial in the higher dimensions.

\smallskip

We need two more  lemmas to take care of some missing cases.

\begin{lemma}
The fundamental representation $V_{\varpi_{2n-1}}$ for $SO_{4n}(\C)$ is good
provided $n \geq 2$.
\end{lemma}

\begin{proof}
We have $\bigwedge^{2n-1}(\C^{4n})|SO_{4n-2}(\C) \times SO_2(\C)$ 
contains the nonself-dual representation  $\bigwedge^{2n-1}(\C^{4n-2})_+
\boxtimes \bigwedge^0(\C^2)$.
\end{proof}

We also need

\begin{lemma}
The fundamental representation $V_{\varpi_{4}}$ for $SO_{4n}(\C)$ is good
provided $n \geq 3$.
\end{lemma}

\begin{proof}
We have $\bigwedge^{4}(\C^{4n})|SO_{4n-2}(\C) \times SO_2(\C)$
contains the representation $\bigwedge^{2}(\C^{4n-2})\boxtimes 
\bigwedge^{2}(\C^{2})$.
\end{proof}

 All other fundamental representations for $SO_{4n}(\C)$ follow from
 Lemma \ref{orthogonalsquarecase}.
 
We have proved the following 

\begin{proposition} \label{specialorthogonal}
All the fundamental representations except the first of $SO_{4n}(\C)$ are good
provided $n\geq 3$.
\end{proposition}

\subsection{The orthogonal groups $SO_{4n+2}(\C)$}
In the case of $SO_{4n+2}(\C)$ the last two fundamental representations
are not self-dual. The Cartan product of the last two fundamental
representations is the exterior power $\bigwedge^{2n}(\C^{4n+2})$.

Every self-dual representation is a Cartan product
of Cartan powers of the first $2n-1$ fundamental representations together with
the Cartan powers of $\bigwedge^{2n}(\C^{4n+2})$.
By Theorem \ref{Cartanproduct} we will be done once we prove
\begin{lemma}
The representation $\bigwedge^{2n}(\C^{4n+2}) = V_{\varpi_{2n}+\varpi_{2n+1}} $ 
is good if $n \geq 2$.
\end{lemma}

\begin{proof}
We have 
$\bigwedge^{2n}(\C^{4n+2})|SO_{4n}(\C) \times SO_2(\C)$
contains $\bigwedge^{2n-2}(\C^{4n}) \boxtimes \bigwedge^{2}(\C^{2})$.
The first factor is good provided $n \geq 2$.
\end{proof}

We have concluded our analysis of the even special orthogonal
groups $SO_{2n}(\C)$.

\begin{proposition}
All of the self-dual representations except the Cartan powers of the
standard representation of $SO_{2n}(\C)$ are good provided $n \geq 5$.
\end{proposition}

\subsection{The exceptional cases for $Sp_{2n}(\C)$ and $SO_{2n+1}(\C)$.}
We first prove that $\overline{\Theta}$ is trivial for torus quotients of 
the space of lines in the symplectic vector space
$\C^{2n}$. Let $x_i,1 \leq i \leq 2n$ be the linear coordinates
relative to an adapted basis chosen as before so $(e_i, e_{2n+1 -i}) =1$
and all other symplectic products are zero. We let
$\overline{i} = 2n+1 -i, 1 \leq i \leq 2n$. It is then apparent
that in any $H$--invariant monomial $x_I$ in the coordinates $x_i$ the indices
$i$ and $\overline{i}$ must appear the same number of times and consequently
$x_I$ is invariant under $\theta$. An analogous argument takes care
of the torus quotients of the quadrics $\mathcal{Q}$. 
We leave to the reader the task of checking that any 
$H$--invariant monomial
in the Pl\"ucker coordinates for $Gr_{2}^{0}(\C^{4})$ (the symplectic case)
and $Gr_{2}^{0}(\C^{5})$ (the orthogonal case) is invariant under $\theta$.

Finally it remains to treat the case of the flag manifold of lines
and planes in $\C^4$. We will do this at the end of 
 \S \ref{generalparameters} for the case of  general {\em real} parameters.

\subsection{The exceptional cases for $SO_{2n}(\C)$.}

We prove that $\overline{\Theta}$ acts trivially on the torus quotients of
the quadrics $\mathcal{Q}$ with a symplectic form
corresponding to an integral orbit in the same
way as we did for torus quotients of the space of lines in the symplectic 
vector space $\C^{2n}$. 

Also $\overline{\Theta}$ acts trivially on the torus quotients of  $Gr_2^0(\C^6)$
because under the isomorphism between $SO_6(\C)$ and $SL_4(\C)$ the
Grassmannian $Gr_2^0(\C^6)$ corresponds to the flag manifold $F_{1,3}(\C^4)$.
Similarly $\overline{\Theta}$ does not act trivially on 
the torus quotient of $F_{1,2}^0(\C^6)$ 
because under the above isomorphism $F_{1,2}^0(\C^6)$ corresponds to
the {\em full} flag manifold $F_{1,2,3}(\C^4)$.

Next we explain why $\overline{\Theta}$ is trivial on the torus
quotients of the Lagrangian spaces 
$Gr_2^0(\C^4)^+$, $Gr_2^0(\C^4)^-$, $Gr_4^0(\C^8)^+$, and $Gr_4^0(\C^8)^-$.
It is easy to check (for example by using the Spin representation) that
$Gr_2^0(\C^4)^+$, and $Gr_2^0(\C^4)^-$ are isomorphic to $\mathbb{CP}^1$
and the corresponding torus quotients are points so it is clear that
$\overline{\Theta}$ is trivial for these two cases. As for
the cases of $Gr_4^0(\C^8)^+$, and $Gr_4^0(\C^8)^-$ it follows from
the triality isomorphism, \cite{FultonHarris}, \S 20.3, that each of these 
two flag manifolds
is isomorphic to the quadric $\mathcal{Q}_6$ by an isomorphism that
is torus equivariant (though perhaps with a different but equivalent action)
and conjugates the Chevalley involution to a new involution that
still acts on the torus by inversion. Hence by Lemma 
\ref{characterizationofChevalley} the new involution is conjugate to the Chevalley involution
by an element $Ad h$ and induces the same map as the Chevalley involution
on any torus quotient.
Since we have seen that $\overline{\Theta}$ is trivial on  $\mathcal{Q}_6//_0 H$
it follows that  $\overline{\Theta}$ is trivial on $Gr_4^0(\C^8)^+//_0 H$
and $Gr_4^0(\C^8)^-//_0 H$.
 \medskip
  
We have now dealt with the cases of $Gr_n^0(\C^{2n})^+$ and
$Gr_n^0(\C^{2n})^-$ for $n=2$ and $n=4$. The torus quotient of
the flag manifold $F_{1,2}^0(\C^4)$ is equal to a point so
it is trivial that $\overline{\Theta}$ is the identity for this case.

\medskip

It remains to prove that $\overline{\Theta}$ is not equal to the identity for 
the quotients$F_{1,2}^0(\C^4)$ and  $F_{1,4}^0(\C^8)^{\pm}//_{0} H$. We will 
prove this in the next section. 

\subsection{From integral parameters to general parameters}
\label{generalparameters}

We first observe that by Proposition \ref{reductionproposition}
it suffices (except for a small number of examples) to promote
the nontriviality results obtained above from integral parameters
to general parameters for {\em Grassmannians}. Here there is no
problem. We use the ``scaling trick''. Namely if we scale the
symplectic form by a real number $c$ (thereby multiplying $a$ and $\br$
by $c$ we do not change the torus quotient and we do not change
$\overline{\Theta}$. To be precise suppose the symplectic form
is induced by embedding $M$ into $\mathfrak{k}^{\ast}$ as the orbit
$\mathcal{O}_{\lambda}$. Let $m_c$ be the automorphism of 
$\mathfrak{k}^{\ast}$ given by multiplication by $c$. Then $m_c$
is $T$--equivariant and we have the following diagram

\begin{center}
\begin{math}
\begin{CD}
\mathcal{O}_{\lambda}//_{\br} T  @>m_c>> \mathcal{O}_{c\lambda}//_{\br} T \\
@A{\overline{\Theta}}AA                  @AA{\overline{\Theta}}A \\
\mathcal{O}_{\lambda}//_{\br} T  @>m_c>> \mathcal{O}_{c\lambda}//_{\br} T.
\end{CD}
\end{math}
\end{center} 
Thus $\overline{\Theta}$ is either trivial for all $c$ or nontrivial for all 
$c$.

This takes care of all the Grassmannian cases
(i.e. the cases where $P$ is maximal). Also by using the 
``scaling trick'' we may promote all of the above results from
integral parameters to rational parameters. By continuity this
allows us to promote all the cases where we have proved that
$\overline{\Theta}$ is trivial to the general case.  

It remains to prove that $\overline{\Theta}$ is not equal to the
identity on the torus quotients of $F_{1,2}^0(\C^4)$ and the
two flag manifolds $F_{1,4}^0(\C^8)^+$ and $F_{1,4}^0(\C^8)^-$. 
In all three cases there is a {\em two} real parameter
family of symplectic forms. Because we have a two parameter
family the scaling trick does not suffice to extend our nontriviality
result from integral parameters to general parameters. We will give complete
details in the second case and third cases
and give the main point for the easier first case. First we claim it
suffices to treat the case of $F_{1,4}^0(\C^8)^+$.  Indeed the matrix in $O(8)$ which
interchanges the fourth and fifth standard basis vectors and leaves
all the other basis vectors fixed interchanges $F_{1,4}^0(\C^8)^+$
and $F_{1,4}^0(\C^8)^-$, commutes with $\Theta$ and normalizes the
torus. Hence $\overline{\Theta}$ is the identity on $F_{1,4}^0(\C^8)^+//_0 T$
if and only if it is the identity on $F_{1,4}^0(\C^8)^-//_0 T$.
We now prove 

\begin{lemma}
The map $\overline{\Theta}$ is not equal to the identity on $F_{1,4}^0(\C^8)^+//_0 T$
for any of the symplectic quotients for the symplectic forms corresponding to 
the orbits of $a \varpi_1 + b \varpi_4=(a+(b/2), b/2,b/2,b/2)$.

\end{lemma}

\begin{proof}

We will apply Lemma \ref{splittingthezerolevel} to deduce that
$\overline{\Theta}$ is not equal to the identity.

\smallskip

Let $St_4^0(\C^8)$ denote the submanifold of the Stiefel manifold of $8$ by $4$ complex
matrices with columns which are orthonormal for the standard
hermitian form  such that the columns span a subspace which is Lagrangian
for the bilinear form 
$(x,y) = \sum_{i=1}^8 x_i y_{9-i}$.  Let 
$\pi_1 : St_4^0(\C^8) \rightarrow Gr_1^0(\C^8)$ be the map sending a matrix to 
the span of its first column, and let 
$\pi_4 : St_4^0(\C^8) \rightarrow Gr_4^0(\C^8)$ take $A$ to $Im(A)$.
Let $\pi_{1,4}: St_4^0(\C^8) \rightarrow F_{1,4}^0(\C^8)$ 
be given by $\pi_{1,4}(A) = (\pi_1(A),\pi_4(A))$.
Let $\omega_1$ be the symplectic form corresponding to $\varpi_1$
for the isotropic Grassmannian $Gr_1^0(\C^8)$, and let $\omega_4$ correspond to 
$2 \varpi_4$ for $Gr_4^0(\C^8)$.  A momentum mapping $\mu_1$ for $Gr_1^0(\C^8)$ is 
given by  
$\mu_1(\pi_1(A)) = (|a_{11}|^2 - |a_{81}|^2, |a_{21}|^2-|a_{71}|^2, |a_{31}|^2 - 
|a_{61}|^2,
|a_{41}|^2-|a_{51}|^2)$.  We claim that 
a momentum mapping for $Gr_4^0(\C^8)$ is given by
$\mu_4(\pi_4(A)) = (|r_1|^2 - |r_8|^2, |r_2|^2-|r_7|^2, |r_3|^2 - |r_6|^2,
|r_4|^2-|r_5|^2)$, where 
$r_i$ is the $i$-th row vector of $A$ and $|r_i|^2$ is the
length of the $i$--th row for the standard Hermitian form on $\C^4$.
Indeed since the columns of $A$ are orthonormal for the standard Hermitian
form on $\C^8$ it follows that the momentum map $\nu$ for the action
of the diagonal torus in $U(8)$ is given by
$\nu(A) = (|r_1|^2, |r_2|^2, \cdots , |r_8|^2).$
The torus $T$ for $SO(8)$ is embedded in the torus for $U(8)$ as the set of
diagonal matrices such
that $z_i = z_{9-i}^{-1}, 1 \leq i \leq 8$. Since the momentum map for $T$
is the orthogonal projection onto the Lie algebra of $T$ the claim follows.   
Thus $\mu_{1,4}^{a,b} = a\mu_1 + b\mu_4$ is a momentum map for the natural 
symplectic embedding 
of $F_{1,4}^0(\C^8)^+$ into $Gr_1^0(\C^8) \times Gr_4^0(\C^8)^+$ for the symplectic 
form $a \omega_1 + b \omega_4$ on the product.

Let 
\[
A = 
\begin{pmatrix}
+\alpha & +\alpha & +\gamma & +\gamma \\
+\alpha & +\alpha & -\gamma & -\gamma \\
+\alpha & -\alpha & +\gamma & -\gamma \\
+\alpha & -\alpha & -\gamma & +\gamma \\
+\beta & +\beta & +\delta & +\delta \\
+\beta & +\beta & -\delta & -\delta \\
-\beta & +\beta & -\delta & +\delta \\
-\beta & +\beta & +\delta & -\delta \\
\end{pmatrix}
\]

First we claim that the flag corresponding to $A$ as above is isotropic
and belongs to $F_{1,4}^0(\C^8)^+$. Indeed the linear functional $x_I$ 
on $\bigwedge^4(\C^8)$ obtained
by wedging together the first four elements of the  basis dual to the
standard basis of $\C^8$ takes value $-16 \alpha^2 \gamma^2$ on $A$
(the determinant of the upper four by four block). But  $x_I$ is fixed
by the Hodge star and consequently takes the value zero on any element
of $F_{1,4}^0(\C^8)^-$.

Next note that the columns of $A \in St_4^0(\C^8)$ 
and $\mu_{1,4}^{a,b}(\pi_{1,4}(A)) = 0$ if 
$\alpha,\beta,\gamma,\delta \in \R$ satisfy the following equations:

\begin{enumerate}
\item $4 \alpha^2 + 4 \beta^2 = 1, 4 \gamma^2 + 4 \delta^2 = 1$ (unit length)
\item $(a+2b)(\alpha^2 - \beta^2) + 2b(\gamma^2 - \delta^2) = 0$ (momentum 0)
\end{enumerate}

For small enough $t \in \R$, solutions are given by  
\begin{align}
\alpha(t) &= \sqrt{1/8 + 2bt} \\ 
\beta(t) &= \sqrt{1/8 - 2bt} \\
\gamma(t) &= \sqrt{1/8 - (a+2b)t} \\
\delta(t) &= \sqrt{1/8 + (a + 2b)t}
\end{align}

For $t \neq 0$, $\alpha(t)^2 \neq \beta(t)^2$, and hence $a \mu_1(\pi_1(A(t))) 
\neq 0$.  Hence, 
$[A(t)] \in F_{1,4}^0(\C^8)^+//_0H$ is not fixed by  
$\overline{\Theta}$ by Lemma \ref{splittingthezerolevel}.

\end{proof}

\vskip 12pt

For the symplectic flag manifold $F_{1,2}^0(\C^4)$, the reader will give an
analogous (but easier) argument using the matrix 

\[
A = 
\begin{pmatrix}
\alpha & \gamma  \\
\alpha & - \gamma  \\
\beta & \delta  \\
\beta & - \delta \\
\end{pmatrix}  
\]

\section{The exceptional groups}

In this section we will prove that the self-duality map $\overline{\Theta}$
is never equal to the identity on a torus quotient of a flag manifold
of an exceptional group.

We will prove the following theorem below using the branching trick

\begin{theorem}
Let $G$ be an exceptional Lie group and $V_{\lambda}$ be a self-dual
representation of $G$. Then $V_{\lambda}$ is good.
\end{theorem}

By Theorem \ref{quantumtoclassical} we then obtain

\begin{corollary}
Let $M$ be an integral self-dual flag manifold associated to an
exceptional group. 
Then the  self-duality map $\overline{\Theta}$ on $M//_0 H$
is not equal to the identity.
\end{corollary}

We then pass from integral parameters to general parameters using
the scaling trick.

\subsection{The group $G_2$}

See \cite{MillsonToledano} Theorem 1.7. By branching to $SL_3(\C)$ 
one finds that both the fundamental representations of $G_2$
are good and hence by Theorem \ref{Cartanproduct} all
representations are good.

\subsection{The group $F_4$}
We restrict the fundamental representations $V_i,1 \leq i \leq 4$ to the
subgroup $Spin(9)$.The only bad representations for $Spin(9)$
are the Cartan powers of the first fundamental representation.
We check by Lie that the restriction of every fundamental representation
contains  a good irreducible summand.  

\subsection{The group $E_6$}
We will use the notation of \cite{Bourbaki}, page 261.

\begin{lemma}
Any self-dual highest weight 
$\lambda$ may be written
$$\lambda = a(\varpi_1 + \varpi_6) + b(\varpi_3 + \varpi_5) + c \varpi_2
+ d \varpi_4.$$
\end{lemma}

As a consequence a self-dual irreducible representation $V_{\lambda}$
is a quotient of a tensor product of Cartan powers
of $V_{\varpi_1 + \varpi_6}$,$V_{\varpi_3 + \varpi_5}$, $V_{\varpi_2}$
and $V_{\varpi_4}$. 

We restrict the four basic representations $V_{\varpi_1 + \varpi_6}$,
$V_{\varpi_3 + \varpi_5}$, $V_{\varpi_2}$ and $V_{\varpi_4}$ to
the maximal rank subgroup $SL_5(\C) \times SL_2(\C)$. We find using LiE that
the restriction of each of the four representations contains either a nonself-dual
or a good irreducible
summand. Hence  each of these representations is good
and hence by
Theorem \ref{Cartanproduct}  any quotient of a  
tensor product involving a Cartan power 
of one of the four basic representations is good. Hence  any self-dual 
representation is good.

\subsection{The group $E_7$}
We restrict the fundamental representations $V_i, 1 \leq i \leq 7$ to
$SL_8(\C)$. The only self-dual representations for $SL_8(\C)$ are the Cartan
powers of $\varpi_1 + \varpi_7$. We again verify by LiE that the restriction of each
$V_i$ all contain either a good representation or a nonself-dual irreducible
summand. 

\subsection{The group $E_8$}
We restrict the fundamental representations to $SL_9(\C)$.
We again check by LiE that
at least one good representation occurs in the restriction of
each fundamental representation of $E_8$.

\section{Further questions}\label{furtherquestions}

It appears that the duality map $\Theta$ preserves almost every important
structure connected with the Grassmannian. We list some that need
to be further investigated.

\subsection{Toric Degenerations of Torus Quotients of Flag Manifolds}

The duality
map carries the toric space that is the degeneration of $G_{m+1}(\C^n)//_{\br}H$
to the toric space that is the degeneration of $G_{n-m-1}(\C^n)//_{\bs}H$
in the construction of P.\ Foth and Yi Hu in \cite{FothHu} since it permutes
the Gelfand-Tsetlin Hamiltonians, see
\S \ref{GTs}. At the moment we have not yet proved that
the duality map can be extended to a map of total spaces of the 
degeneration. There are a number of other toric degenerations of flag manifolds,
which induce toric degenerations of torus quotients of  flag manifolds, see
for example \cite{GonciuleaLakshmibai}, Chapter 11, where it remains to extend 
the
duality map to the total space of the degenerations.

\subsection{Gel'fand Hypergeometric Functions}
In the 1980's Gel'fand and his collaborators created a theory of
hypergeometric functions on Grassmannians generalizing the classical theory 
defined on the moduli spaces $\mathcal{M}_r(\mathbb{CP}^1)$, see
\cite{DeligneMostow}. There is evidence that our duality
is compatible with these functions. Indeed in the early real version of the
theory this is proved in \cite{GelfandGraev}. However,
it is not easy to see how duality of hypergeometric functions would go
in the complex case. Nevertheless, there are reasons to believe that there
should be such a duality. To begin with, 
the duality map preserves the GGMS stratifications, see \cite{GGMS}.
We recall the definition. The GGMS strata are parametrized by 
matroids on the set $1,2,...,n$. 
Two points $x$ and $y$ are
in the same stratum if $M(x) = M(y)$. For $M$ a matroid on $1,2,...,n$ we let
$S_M$ denote the corresponding stratum. The dual $M^{\ast}$ of a
matroid $M$ is defined and discussed in \cite{Oxley}, Chapter 2.
We have the following theorem \cite{HowardMillson}

\begin{theorem}
\hfill
$$\Psi(S_M) = S_{M^{\ast}}.$$
\end{theorem} 
However it is not clear that the duality map lifts to the families of 
arrangement
complements over the strata that give rise to the hypergeometric functions.

If a duality of hypergeometric functions could be established it would
afford the opportunity to carry over the very detailed information on
monodromy obtained in \cite{DeligneMostow} to some of the Gelfand examples.

\subsection{Self-dual torus orbits}
A very concrete problem suggested by the work of \cite{DO} is the 
problem of finding
the fixed points of the duality  map (``self-associated point sets''
in the terminology of \cite{DO}). 
 This problem is discussed in  detail in Chapter III of 
\cite{DO}. We note that in \cite{Foth}, P.\ Foth gave a 
description of the fixed-point set of an { \em anti-holomorphic} involution
on a weight variety. The problem of finding the self-dual torus orbits
of flags is the analogous problem for the (holomorphic) Chevalley involution.
However the results of \cite{DO} suggest this problem will be more subtle.
For example, it is proved in \cite{DO} following \cite{Coble} that in two cases
the fixed set is closely related to the moduli theory of curves.
Since one of
these results is very easy to describe we conclude with it. Note first that we 
obtain
a self-dual torus quotient by giving $Gr_n(\C^{2n})$ the symplectic form
corresponding to $2 \varpi_n$ and taking $\br = (1,1,\cdots,1)=\varpi_{2n}$
so $a=2 = |\br|/n$.
In \cite{DO}, Theorem 4, pg. 51, it is proved that the fixed set $S_{n-1}$ of
the self-duality $\overline{\Theta}$ acting on the resulting torus quotient 
(equivalently the moduli space of  $2n$--tuples
of equally weighted (by $1$) points in $\mathbb{CP}^{n-1}$) is a rational subvariety of dimension
$\frac{n(n+1)}{2}$. In Example 4, pg. 37, of \cite{DO} it is proved  that
$S_2$ is isomorphic to the Baily-Borel-Satake
compactification of the (level two) Siegel modular variety of genus $2$.

\end{document}